\documentclass[a4paper,thmsa]{article}%
\usepackage{graphicx}
\usepackage{amsmath}
\usepackage{amsfonts}
\usepackage{amssymb}
\providecommand{\U}[1]{\protect\rule{.1in}{.1in}}
\newtheorem{theorem}{Theorem}

\newtheorem{lemma}{Lemma}

\newtheorem{proposition}[theorem]{Proposition}

\newenvironment{proof}[1][Proof]{\textbf{#1.} }{\ \rule{0.5em}{0.5em}}
\numberwithin{equation}{section}

\def\de{\overset{\mathrm{def}}{=}}
\def\sni{\smallskip\noindent}

\def\mi\mathfrak{i}

\begin{document}

\title{Lingering random walks in random environment on a strip}
\author{Erwin Bolthausen$^{1}$ and Ilya Goldsheid$^{2}$\\
Universit\"{a}t Z\"{u}rich$^{1}$\\Queen Mary, University of London$^{2}$}
\maketitle

\begin{abstract}
We consider a recurrent random walk (RW) in random environment (RE) on a
strip. We prove that if the RE is i. i. d. and its distribution is not
supported by an algebraic subsurface in the space of parameters defining the
RE then the RW exhibits the $(\log t)^{2}$ asymptotic behaviour. The
exceptional algebraic subsurface is described by an explicit system of
algebraic equations.

One-dimensional walks with bounded jumps in a RE are treated as a
particular case of the strip model. If the one dimensional RE is
i. i. d., then our approach leads to a complete and constructive
classification of possible types of asymptotic behaviour of
recurrent random walks. Namely, the RW exhibits the $(\log t)^{2}$
asymptotic behaviour if the distribution of the RE is not
supported by a hyperplane in the space of parameters which shall
be explicitly described. And if the support of the RE belongs to
this hyperplane then the corresponding RW is a martingale and its
asymptotic behaviour is governed by the Central Limit Theorem.


\end{abstract}

\section{Introduction}

The aim of this work is to describe conditions under which a recurrent random
walk in a random environment (RWRE) on a strip exhibits the $\log^{2} t$
asymptotic behaviour. This slow, lingering movement of a walk was discovered
by Sinai in 1982 \cite{S}. At the time, this work had brought to a logical
conclusion the study of the so called simple RWs (SRW) started by Solomon in
\cite{So} and by Kesten, Kozlov, and Spitzer in \cite{KKS}. The somewhat
misleading term ``simple" is often used as an abbreviation describing a walk
on a one-dimensional lattice with jumps to nearest neighbours.

Our work was motivated by a question asked by Sinai in \cite{S} about the
validity of his (and related) results for other models. Perhaps the simplest
extension of the SRW is presented by a class of one-dimensional walks whose
jumps (say) to the left are bounded and to the right are of length at most
one. These models were successfully studied by a number of authors and the
relevant references can be found in \cite{Br}. We would like to quote one
result concerning this special case since it is perhaps most close to our
results stated below in Theorems \ref{Theor2} and \ref{Theor3}. Namely,
Bremont proved in \cite{Br1} that if the environment is defined by a Gibbs
measure on a sub-shift of finite type, then the asymptotic behaviour of a
recurrent RW is either as in the Sinai's theorem, or it is governed by the
Central Limit Law.

General 1DWBJ were also studied by different authors. Key in \cite{Ke} found
conditions for recurrence of a wide class of 1DWBJ. Certain sufficient
conditions for the Sinai behaviour of 1DWBJ were obtained by Letchikov in
\cite{L}. The results from \cite{L} will be discussed in a more detailed way
in Section \ref{sec1.1} after the precise definition of the one-dimensional
model is given. We refer the reader to \cite{Z} for further historical
comments as well as for a review of other recent developments.

The main object of this paper is the RWRE on a strip. We prove (and this is
the main result of this paper) that recurrent walks in independent identically
distributed (i. i. d.) random environments on a strip exhibit the $\log^{2} t$
asymptotic behaviour if the support of the distribution of the parameters
defining the random environment does not belong to a certain algebraic
subsurface in the space of parameters. This subsurface is defined by an
explicit system of algebraic equations.

The one dimensional RW with bounded jumps can be viewed as a
particular case of a RWRE on a strip. This fact was explained in
\cite{BG} and we shall repeat this explanation here. Due to this
reduction, our main result implies a complete classification of
recurrent 1DWBJ in i.i.d. environments. Namely, the corresponding
system of algebraic equations reduces in this case to one linear
equation which defines a hyperplane in the space of parameters. If
the support of the distribution of parameters does not belong to the
this hyperplane, then the RW exhibits the Sinai behaviour (see
Theorem \ref{Theor2} below). But if it does, then (Theorem
\ref{Theor3} below) the corresponding random walk is a martingale
and its asymptotic behaviour is governed by the Central Limit Law.
In brief, recurrent 1DWBJ are either of the Sinai type, or they are
martingales.

In the case of a strip, a complete classification can also be obtained and it
turns out that once again the asymptotic behaviour is either the Sinai, or is
governed by the Invariance Principle. However, this case is less transparent
and more technical even to describe in exact terms and we shall leave it for a
future work.

\bigskip\noindent The paper is organized as follows. We state Sinai's result
and define a more general one-dimensional model in Section
\ref{sec1.1}. Section \ref{sec1.2} contains the definition of the
strip model and the explanation of the reduction of the
one-dimensional model to the strip case. Main results are stated in
Section \ref{mainresults}. Section \ref{sec2} contains several
statements which are then used in the proof of the main result,
Theorem \ref{Theor1}. In particular, we introduce random
transformations associated with random environments in Section
\ref{randomtransformations}. It turns out to be natural to recall
and to extend slightly, in the same Section
\ref{randomtransformations}, those results from \cite{BG} which are
used in this paper. An important Lemma \ref{alg-equiv} is proved in
Section \ref{sec2.3}; this Lemma allows us to present the main
algebraic statement of this work in a constructive form. In section
\ref{sec2.4} we prove the invariance principle for the $\log$ of a
norm of a product of certain matrices. This function plays the role
of the so called potential of the environment and is responsible for
the Sinai behaviour of the random walk. It is used in the proof of
our main result in Section \ref{sec3.1}.

Finally Appendix contains results of which many (if not all) are not
new but it is convenient to have them in a form directly suited for
our purposes. Among these, the most important for our applications
is the Invariance Principle (IP) for ``contracting" Markov chains
(Section \ref{rand-tr}). Its proof is derived from a well known IP
for general Markov chains which, in turn, is based on the IP for
martingales.

\textit{Conventions.} The following notations and terminology shall be used
throughout the paper. $\mathbb{R}$ is the set of real numbers, $\mathbb{Z}$ is
the set of integer numbers, and $\mathbb{N}$ is the set of positive integers.

For a vector ${x}=(x_{i})$ and a matrix $A=(a(i,j))$ we put
\[
\left\|  {x}\right\|  \overset{\mathrm{def}}{=}\max_{i}|x_{i}|,\ \left\|
A\right\|  \overset{\mathrm{def}}{=}\max_{i}\sum_{j}|a(i,j)|.
\]
Note that $\left\|  A\right\|  = \sup_{\left|  {x}\right|  =1}\left\|
A{x}\right\|  $. We say that $A$ is strictly positive (and write $A>0),$ if
all its matrix elements satisfy $a(i,j)>0$. $A$ is called non-negative (and we
write $A\geq0)$, if all $a(i,j)$ are non negative. A similar convention
applies to vectors.

\subsection{Sinai's result and some of its extensions to 1DWBJ. \label{sec1.1}%
}

Let $\omega\overset{\mathrm{def}}{=}(p_{n})_{-\infty<n<\infty}$ be a sequence
of independent identically distributed (i. i. d.) random variables, satisfying
$\varepsilon\le p_{n}\le1-\varepsilon$, where $\varepsilon>0$. Put $q_{n}=
1-p_{n}$ and consider a random walk $\xi(t)$ on a one-dimensional lattice with
a starting point $\xi(0)=0$ and transition probabilities
\[
Pr_{\omega}\{\,\xi(t+1)=n+1\,|\, \xi(t)=n\,\}= p_{n},\ \ Pr_{\omega}%
\{\,\xi(t+1)=n-1\,|\, \xi(t)=n\,\}= q_{n}
\]
thus defining a measure $Pr_{\omega}\{\cdot\}$ on the space of trajectories of
the walk. It is well known (Solomon, \cite{So}) that this RW is recurrent in
almost all environments $\omega$ if and only if $\mathbb{E}\ln\frac{q_{n}%
}{p_{n}}=0$ (here $\mathbb{E}$ denotes the expectation with respect to the
relevant measure $\mathbb{P}$ on the space of sequences). In \cite{S} Sinai
proved that if $\mathbb{E}(\ln\frac{q_{n}}{p_{n}})^{2}>0$ and $\xi(\cdot)$ is
recurrent then there is a weakly converging sequence of random variables
$b_{t}(\omega),\ t=1,2,...$ such that
\begin{equation}
\label{0.1}(\log t)^{-2}\xi(t)-b_{t}{\to} 0 \ \hbox{ as }\ \ t\to\infty.
\end{equation}
The convergence in (\ref{0.1}) is in probability with respect to the so called
annealed probability measure $\mathbb{P}(d\omega)Pr_{\omega}$ (for precise
statements see section \ref{mainresults}). The limiting distribution of
$b_{t}$ was later found, independently, by Golosov \cite{G1,G2} and Kesten
\cite{K}.

The one-dimensional walk with bounded jumps on $\mathbb{Z}$ is defined
similarly to the simple RW. Namely let ${\omega} \overset{\mathrm{def}}{=}
(p(n,\cdot))$, $n\in\mathbb{Z}$, be a sequence of non-negative vectors with
$\sum_{k=-m}^{m}p(n,k)=1$ and $m>1$. Put $\xi(0)=0$ and
\begin{equation}
\label{1dtransition}Pr_{{\omega}}\left(
\xi(t+1)=n+k\,|\,\xi(t)=n\right)
\overset{\mathrm{def}}{=}p(n,k),\quad n\in\mathbb{Z}.
\end{equation}
Suppose next that $p(n,\cdot)$ is a random stationary in $n$ (in particular it
can be i. i. d.) sequence of vectors. Sinai's question can be put as follows:
given that a RW is recurrent, what kind of asymptotic behaviour would one
observe? and under what conditions?

There were several attempts to extend Sinai's result to the
(\ref{1dtransition}) model. In particular, Letchikov \cite{L} proved that if
for some $\varepsilon>0$ with $\mathbb{P}$-probability 1
\[
p(n,1)\ge\sum_{k=-m}^{-2}p(n,k)+\varepsilon\ \hbox{ and }\ p(n,-1)\ge
\sum_{k=2}^{m}p(n,k)+\varepsilon
\]
and the distribution of the i. i. d. random vectors $p(n,\cdot)$ is absolutely
continuous with respect to the Lebesgue measure (on the relevant simplex),
then the analogue of Sinai's theorem holds. (In \cite{L}, there are also other
restrictions on the distribution of the RE but they are much less important
than the ones listed above.)

The technique we use in this work is completely different from that used in
\cite{Ke}, \cite{L}, \cite{Br}, \cite{Br1}. It is based on the methods from
\cite{BG} and \cite{G} and this work presents further development of the
approach to the analysis of the RWRE on a strip started there.

\subsection{Definition of the strip model. \label{sec1.2}}

The description of the strip model presented here is the same as in \cite{BG}.

Let $(P_{n},Q_{n},R_{n}),\ -\infty<n<\infty,$ be a strictly stationary ergodic
sequence of triples of $m\times m$ matrices with non-negative elements such
that for all $n\in\mathbb{Z}$ the sum $P_{n}+Q_{n}+R_{n}$ is a stochastic
matrix,
\begin{equation}
(P_{n}+Q_{n}+R_{n})\mathbf{1}=\mathbf{1}, \label{stoch}%
\end{equation}
where $\mathbf{1}$ is a column vector whose components are all equal to $1$.
We write the components of $P_{n}$ as $P_{n}(i,j),$ $1\leq i,j\leq m,$ and
similarly for $Q_{n}$ and $R_{n}.$ Let $(\Omega,\mathcal{F},\mathbb{P}%
,\mathcal{T})$ be the corresponding dynamical system with $\Omega$ denoting
the space of all sequences $\omega=(\omega_{n})=((P_{n},Q_{n},R_{n}))$ of
triples described above, $\mathcal{F}$ being the corresponding natural
$\sigma$-algebra, $\mathbb{P}$ denoting the probability measure on
$(\Omega,\mathcal{F})$, and $\mathcal{T}$ being a shift operator on $\Omega$
defined by $(\mathcal{T}\omega)_{n}=\omega_{n+1}$. For fixed $\omega$ we
define a random walk $\xi(t),$ $t\in\mathbb{N}$ on the strip $\mathbb{S}%
=\mathbb{Z}\times\{1,\ldots,m\}$ by its transition probabilities
$\mathcal{Q}_{\omega}(z,z_{1})$ given by
\begin{equation}
\mathcal{Q}_{\omega}(z,z_{1})\overset{\mathrm{def}}{=}\left\{
\begin{array}
[c]{ll}%
P_{n}(i,j) & \mathrm{if\quad}z=(n,i),\ z_{1}=(n+1,j),\\
R_{n}(i,j) & \mathrm{if\quad}z=(n,i),\ z_{1}=(n,j),\\
Q_{n}(i,j) & \mathrm{if\quad}z=(n,i),\ z_{1}=(n-1,j),\\
0 & \mathrm{otherwise,}%
\end{array}
\right.  \label{striptransition}%
\end{equation}
This defines, for any starting point $z=(n,i)\in\mathbb{S}$ and any $\omega$,
a law $Pr_{\omega,z}$ for the Markov chain $\xi(\cdot)$ by
\begin{equation}
Pr_{\omega,z}\left(  \xi(1)=z_{1},\ldots,\xi(t)=z_{t}\right)  \overset
{\mathrm{def}}{=}\mathcal{Q}_{\omega}(z,z_{1})\mathcal{Q}_{\omega}(z_{1},
z_{2})\cdots\mathcal{Q}_{\omega}(z_{t-1},z_{t}). \label{StripRWRE}%
\end{equation}

We call $\omega$ the \textit{environment }or the \textit{random environment}
on a strip $\mathbb{S}$. Denote by $\Xi_{z}$ the set of trajectories
$\xi(\cdot)$ starting at $z$. $Pr_{\omega,z}$ is the so called quenched
probability measure on $\Xi_{z}$. The semi-direct product $\mathbb{P}%
(d\omega)Pr_{\omega,z}(d\xi)$ of $\mathbb{P}$ and $Pr_{\omega,z}$ is defined
on the direct product $\Omega\times\Xi_{z}$ and is called the annealed
measure. All our main results do not depend on the choice of the starting
point $z$. We therefore write $Pr_{\omega}$ instead of $Pr_{\omega,z}$ when
there is no danger of confusion.


\smallskip\noindent The one-dimensional model (\ref{1dtransition}) reduces to
a RW on a strip due to the following geometric construction. Note
first that it is natural to assume (and we shall do so) that at
least one of the following inequalities holds:
\begin{equation}
\label{1.6}\mathbb{P}\{{\omega}\,:\,p(x,m)>0\}>0\ \hbox{ or }\ \mathbb{P}%
\{{\omega}\,:\,p(x,-m)>0\}>0.
\end{equation}
Consider the one-dimensional lattice as a subset of the $X$-axis in a
two-dimensional plane. Cut this axes into equal intervals of length $m$ so
that each of them contains exactly $m$ consecutive integer points. Turn each
of these intervals around its left most integer point anti-clockwise by
$\pi/2$. The image of $\mathbb{Z}$ obtained in this way is a part of a strip
with distances between layers equal to $m$. Re-scaling the $X$-axis of the
plane by $m^{-1}$ makes the distance between the layers equal to one. The
random walk on the line is thus transformed into a random walk on a strip with
jumps to nearest layers.

The formulae for matrix elements of the corresponding matrices $P_{n}%
,Q_{n},R_{n}$ result now from a formal description of this construction.
Namely, present $x\in\mathbb{Z}$ as $x=nm+i$, where $1\leq i\leq m$. This
defines a bijection $x\leftrightarrow(n,i)$ between the one-dimensional
lattice $\mathbb{Z}$ and the strip $\mathbb{S}=\mathbb{Z}\times\{1,\ldots
,m\}.$ This bijection naturally transforms the $\xi$-process on $\mathbb{Z}$
into a walk on $\mathbb{Z}\times\{1,\ldots,m\}$. The latter is clearly a
random walk of type (\ref{StripRWRE}) and the corresponding matrix elements
are given by
\begin{equation}
\label{1dmodel}\begin{aligned} P_{n}(i,j)=&p(nm+i,m+j-i),\\ R_{n}(i,j)=&p(nm+i,j-i),\\ Q_{n}(i,j)=&p(nm+i,-m+j-i). \end{aligned}
\end{equation}

\subsection{Main results. \label{mainresults}}

Denote by $\mathcal{J}$ the following set of triples of $m\times m$ matrices:
\[
\mathcal{J}\overset{\mathrm{def}}{=}\left\{  (P,Q,R)\,:\,P\ge0,\,
Q\ge 0,\,R\ge0 \ \hbox{ and }\
(P+Q+R)\mathbf{1}=\mathbf{1}\right\}.
\]
Let $\mathcal{J}_{0}\subset\mathcal{J}$ be the support of the probability
distribution of the random triple $(P_{n},Q_{n},R_{n})$ defined above
(obviously, this support does not depend on $n$). The two assumptions
$\mathbf{C1}$ and $\mathbf{C2}$ listed below will be referred to as Condition
$\mathbf{C}$.

\begin{description}
\item[Condition $\mathbf{C}$]

\item[C1] $(P_{n},Q_{n},R_{n})$, $-\infty<n<\infty$, is a sequence of
independent identically distributed random variables.

\item[$\mathbf{C2}$] There is an $\varepsilon>0$ and a positive integer number
$l<\infty$ such that for any $(P,Q,R)\in\mathcal{J}_{0}$ and all
$i,\,j\in[1,m]$
\[
||R^{l}||\le1-\varepsilon,\ \ ((I-R)^{-1}P)(i,j)\ge\varepsilon,
\ \ ((I-R)^{-1}Q)(i,j)\ge\varepsilon.
\]

\end{description}

\smallskip\noindent\textit{Remarks.} 1. We note that say
$((I-R_{n})^{-1} P_{n})(i,j)$ is the probability for a RW starting
from $(n,i)$ to reach $(n+1,j)$ at its first exit from layer $n$.
The inequality $||R_{n}^{l}||\le 1-\varepsilon$ is satisfied in
essentially all interesting cases and, roughly speaking, means
that the probability for a random walk to remain in layer $n$
after a certain time $l$ is small uniformly with respect to $n$
and $\omega$.

2. If the strip model is obtained from the one-dimensional model, then
$\mathbf{C2}$ may not be satisfied by matrices (\ref{1dmodel}). This
difficulty can be overcome if we replace $\mathbf{C2}$ by a much milder
condition, namely:

\begin{description}
\item[$\mathbf{C3}$] For $\mathbb{P}$ - almost all $\omega$:

(a) the strip $\mathbb{S}$ is the (only) communication class of the walk,

(b) there is an $\varepsilon>0$ and a triple $(P,Q,R)\in\mathcal{J}_{0}$ such
that at least one of the following two inequalities holds: $((I-R)^{-1}%
P)(i,j)\ge\varepsilon$ for all $i,\,j\in[1,m]$, or $((I-R)^{-1}Q)(i,j)\ge
\varepsilon$ for all $i,\,j\in[1,m]$.
\end{description}

Our proofs will be carried out under Condition $\mathbf{C2}$. They
can be modified so that to make them work also under Condition
$\mathbf{C3}$. Lemma \ref{contraction} which is used in the proof of
Theorem \ref{Theor1} is the main statement requiring a more careful
treatment under condition $\mathbf{C3}$ and the corresponding
adjustments are not difficult. However, the proofs become more
technical in this case, and we shall not do this in the present
paper. If now vectors $p(x,\cdot)$ defining matrices (\ref{1dmodel})
are $\mathbb{P}$-almost surely such that $p(x,1)\ge\epsilon$ and
$p(x,-1)\ge\epsilon$ for some $\epsilon>0$, then it is easy to see
that Condition $\mathbf{C3}$ is satisfied. We note also that if in
addition the inequalities $p(x,m)\ge\epsilon$ and
$p(x,-m)\ge\epsilon$ hold $\mathbb{P}$-almost surely, then also
$\mathbf{C2}$ is satisfied.

For a triple of matrices $(P,Q,R)\in\mathcal{J}_{0}$ denote by $\pi
=\pi_{(P,Q,R)}=(\pi_{1},\ldots,\pi_{m})$ a row vector with non-negative
components such that
\[
\pi(P+Q+R)=\pi\ \hbox{ and } \sum_{j=1}^{m}\pi_{j}=1.
\]
Note that the vector $\pi$ is uniquely defined. Indeed, the equation for $\pi$
can be rewritten as
\[
\pi(I-R)\left(  (I-R)^{-1}P+(I-R)^{-1}Q\right)  =\pi(I-R).
\]
According to condition $\mathbf{C2}$, the stochastic matrix $(I-R)^{-1}%
P+(I-R)^{-1}Q$ has strictly positive elements (in fact they are $\ge
2\varepsilon$). Hence $\pi(I-R)$ is uniquely (up to a multiplication by a
number) defined by the last equation and this implies the uniqueness of $\pi$.

Consider the following subset of $\mathcal{J}$:
\begin{equation}
\label{algebra}\mathcal{J}_{al}\overset{\mathrm{def}}{=}\{\,(P,Q,R)\in
\mathcal{J}\,:\, \pi(P-Q)\mathbf{1}=0, \ \hbox{where}\ \ \pi(P+Q+R)=\pi\ \,\},
\end{equation}
where obviously $\pi(P-Q)\mathbf{1}\equiv\sum_{i=1}^{m}\pi_{i}\sum_{j=1}%
^{m}(P(i,j)-Q(i,j)). $ Note that $\mathcal{J}_{al}$ is an algebraic subsurface
in $\mathcal{J}$.

We are now in a position to state the main result of this work:

\begin{theorem}
\label{Theor1} Suppose that Condition $\mathbf{C}$ is satisfied, the random
walk $\xi(\cdot)=(X(\cdot),Y(\cdot))$ is recurrent, and $\mathcal{J}%
_{0}\not \subset \mathcal{J}_{al}$. Then there is a sequence of random
variables $b_{t}(\omega),\ t=1,2,...$, which converges weakly as $t\to\infty$
and such that for any $\epsilon>0$
\begin{equation}
\label{sinai}\mathbb{P}\left\{  \omega\,:\, Pr_{\omega}\left(  |\frac
{X(t)}{(\log t)^{2}}-b_{t}|\le\epsilon\right)  \ge1-\epsilon\right\}  \to1
\hbox{ as
}\ t\to\infty.
\end{equation}

\end{theorem}

\textit{Remark.} The algebraic condition in this Theorem requires a certain
degree of non-degeneracy of the support $\mathcal{J}_{0}$ of the distribution
of $(P_{n},Q_{n},R_{n})$. It may happen that relations (\ref{sinai}) hold even
when $\mathcal{J}_{0}\subset\mathcal{J}_{al}$. However Theorem \ref{Theor3}
shows that there are important classes of environments where relations
(\ref{sinai}) (or (\ref{sinai1})) hold if and only if this non-degeneracy
condition is satisfied.

We now turn to the one-dimensional model. It should be mentioned right away
that Theorem \ref{Theor2} is essentially a corollary of Theorem \ref{Theor1}.

Denote by $\tilde{\mathcal{J}}$ the set of all $2m+1$-dimensional probability
vectors:
\[
\tilde{\mathcal{J}}\overset{\mathrm{def}}{=}\{(p(j))_{-m\le j\le m}\,: \,
p(\cdot)\ge0 \ \hbox{ and }\ \sum_{j=-m}^{m} p(j)=1\ \}.
\]
Remember that in this model the environment is a sequence of vectors:
$\omega=\left(  p(x,\cdot)\right)  _{-\infty<x<\infty}$, where $p(x,\cdot
)\in\tilde{\mathcal{J}}$. Let $\tilde{\mathcal{J}}_{0}\subset\tilde
{\mathcal{J}}$ be the support of the distribution of the random vector
$p(0,\cdot)$. Finally, put
\begin{equation}
\label{algebra1}\tilde{\mathcal{J}}_{al}\overset{\mathrm{def}}{=}%
\{\,p(\cdot)\in\tilde{\mathcal{J}}\,:\, \ \sum_{j=-m}^{m}jp(j)=0\ \,\}.
\end{equation}
\begin{theorem}
\label{Theor2} Suppose that:

\noindent(a) $p(x,\cdot),\ x\in\mathbb{Z}$, is a sequence of i. i. d. vectors,

\noindent(b) there is an $\varepsilon>0$ such that
$p(0,1)\ge\varepsilon$, $p(0,-1)\ge\varepsilon$,
$p(0,m)\ge\varepsilon$, and $p(0,-m)\ge\varepsilon$ for any
$p(0,\cdot)\in\tilde{\mathcal{J}}_{0}$,

\noindent(c) for $\mathbb{P}$ almost all environments $\omega$ the
corresponding one-dimensional random walk $\xi(\cdot)$ is recurrent,

\noindent(d) $\tilde{\mathcal{J}}_{0}\not \subset \tilde{\mathcal{J}}_{al}$.

Then there is a weakly converging sequence of random variables $b_{t}%
(\omega),\ t=1,2,...$ such that for any $\epsilon>0$
\begin{equation}
\label{sinai1}\mathbb{P}\left\{  \omega\,:\, Pr_{\omega}\left(  |\frac{\xi
(t)}{(\log t)^{2}}-b_{t}|\le\epsilon\right)  \ge1-\epsilon\right\}  \to1
\hbox{ as
}\ t\to\infty.
\end{equation}

\end{theorem}

\textbf{Proof.} Since the one-dimensional model reduces to a model
on a strip, the result in question would follow if we could check
that all conditions of Theorem \ref{Theor1} follow from those of
Theorem \ref{Theor2}.

It is obvious from formulae (\ref{1dmodel}) that the i. i. d.
requirement (Condition $\mathbf{C1}$) follows from condition $(a)$
of Theorem \ref{Theor2}. We have already mention above that and
Condition $\mathbf{C2}$ follows from condition $(b)$. The recurrence
of the corresponding walk on a strip is also obvious.

Finally, condition $(d)$ implies the algebraic condition of Theorem
\ref{Theor1}. Indeed, formulae (\ref{1dmodel}) show that matrices
$P_{n}$, $Q_{n}$, $R_{n}$ are defined by probability vectors
$p(nm+i,\cdot )\in\tilde{\mathcal{J}}_{0}$, where $1\le i\le m$. Put
$n=0$ and choose all these vectors to be equal to each other, say
$p(i,\cdot)=p(\cdot)\in \tilde{\mathcal{J}}_{0}$, where $1\le i\le
m$. A direct check shows that the triple of matrices $(P,Q,R)$ built
from this vector has the property that $P+Q+R$ is double-stochastic
and irreducible (irreducibility follows from the conditions
$p(1)\ge\varepsilon$ and $p(-1)\ge\varepsilon$). Hence the only
probability vector $\pi$ satisfying $\pi(P+Q+R)=\pi$ is given by
$\pi =(m^{-1},...,m^{-1})$. One more direct calculation shows that
in this case
\[
m\pi(P-Q)\mathbf{1}=\sum_{j=-m}^{m}jp(j).
\]
Hence the condition $\mathcal{J}_{0}\not \subset \mathcal{J}_{al}$ of Theorem
\ref{Theor1} is satisfied if there is at least one vector $p(\cdot)\in
\tilde{\mathcal{J}}_{0}$ such that $\sum_{j=-m}^{m}jp(j)\not =0$. $\Box$




We conclude this section with a theorem which shows, among other
things, that the algebraic condition of Theorem \ref{Theor2} is also
necessary for having (\ref{sinai1}). This theorem does not require
independence as such but in a natural sense it finalizes the
classification of the one-dimensional recurrent RWs with bounded
jumps in the i. i. d. environments.

\begin{theorem}
\label{Theor3} Consider a one-dimensional RW and suppose that

\noindent(a) $p(x,\cdot),\ x\in\mathbb{Z}$, is a strictly stationary ergodic
sequence of vectors,

\noindent(b) there is an $\varepsilon>0$ such that $p(0,1)\ge\varepsilon$ and
$p(0,-1)\ge\varepsilon$ for any $p(0,\cdot)\in\tilde{\mathcal{J}}_{0}$,

\noindent(c) $\tilde{\mathcal{J}}_{0}\subset\tilde{\mathcal{J}}_{al}$, that
is
\[
\sum_{j=-m}^{m} jp(j)=0\ \hbox{ for any $p(\cdot)\in
\tilde{\mathcal{J}}_0$ }.
\]
Then:

\noindent(i) The random walk $\xi(\cdot)$ is asymptotically normal in every(!)
environment $\omega=\left(  p(x,\cdot)\right)  _{\ -\infty< x< \infty}$.

\noindent(ii) There is a $\sigma>0$ such that for $\mathbb{P}$-a. e. $\omega$
\begin{equation}
\label{clt}\lim_{t\to\infty}Pr_{\omega}\left\{  \frac{\xi(t)}{\sqrt{t}}\le
x\right\}  =\frac{1}{\sqrt{2\pi}\sigma}\int_{-\infty}^{x} e^{-\frac{u^{2}%
}{2\sigma^{2}}}du,
\end{equation}
where $x$ is any real number and the convergence in (\ref{clt}) is uniform in
$x$.
\end{theorem}

\textbf{Remarks about the proof of Theorem \ref{Theor3}.} The condition of
this Theorem implies that $\xi(t)$ is a martingale:
\[
E_{\omega}(\xi(t)-\xi(t-1)\,|\,\xi(t-1)=k) =\sum_{j=-m}%
^{m}jp(k,j)=0,
\]
where $E_{\omega}$ denotes the expectation with respect to the
probability measure $Pr_{\omega}$ on the space of trajectories of
the random walk (we assume that $\xi(0)=0$). Let
$U_{n}=\xi(n)-\xi(n-1)$ and put
\[
\sigma_{n}^{2}\overset{\mathrm{def}}{=}E_{\omega}(U_{n}^{2}%
\,|\,\xi(n-1)) =\sum_{j=-m}^{m}j^{2}p(\xi(n-1),j).
\]
Obviously $\varepsilon\le\sigma_{n}^{2}\le m^{2}$, where $\varepsilon$ is the
same as in Theorem \ref{Theor3}. Next put $V_{n}^{2}\overset{\mathrm{def}}%
{=}\sum_{j=1}^{n}\sigma_{j}^{2}$ and
$s_{n}^{2}\overset{\mathrm{def}}{=}
E_{\omega}(V_{n}^{2})={E}_{\omega}(\xi(n)^{2})$. It is useful to
note that $n\varepsilon\le V_{n}^{2},\ s_{n}^{2}\le nm^{2}$. Let $T_{t}%
=\inf\{n:V_{n}^{2}\ge t\}$.

Statement $(i)$ of Theorem \ref{Theor3} is a particular case of a much more
general theorem of Drogin who in particular proves that $t^{-1/2} \xi(T_{t})$
converges weakly to a standard normal random variable. We refer to \cite{HH},
page 98 for more detailed explanations.

Statement $(ii)$ of Theorem \ref{Theor3} is similar to a well known result by
Lawler \cite{La}. The main ingredient needed for proving $(ii)$ is the
following claim:
\begin{equation}
\label{1.13}\hbox{The limit } \lim_{n\to\infty}n^{-1}V_{n}^{2} =\lim
_{n\to\infty}n^{-1} s_{n}^{2}\ \hbox{ exist for $\mathbb{P}$-almost
all }\ \omega.
\end{equation}
Once this property of the variance of $\xi(\cdot)$ is established, $(ii)$
becomes a corollary of Brown's theorem (see Theorems \ref{CLT} and \ref{IP} in
Appendix or Theorem 4.1 in \cite{HH}).

However proving (\ref{1.13}) is not an entirely straightforward matter. The
proof we are aware of uses the approach known under the name ``environment
viewed from the particle". This approach was used in \cite{La} for proving
properties of variances similar to (\ref{1.13}); unfortunately, the conditions
used in \cite{La}, formally speaking, are not satisfied in our case.
Fortunately, Zeitouni in \cite{Z} found the way in which Lawler's result can
be extended to more general martingale-type random walks in random
environments which include our case. $\Box$

\section{Preparatory results.\label{sec2}}

\subsection{Elementary corollaries of condition $\mathbf{C}$.\label{sec2.1}}

We start with several elementary observations following from $\mathbf{C2}$.
Lemma \ref{lemma3} and a stronger version of Lemma \ref{lemma1} can be found
in \cite{BG}. Lemmas \ref{lemma2} and \ref{lemma4} are borrowed from \cite{G}.

\begin{lemma}
\label{lemma1} If Condition $\mathbf{C2}$ is satisfied then for $\mathbb{P}%
$-almost every environment $\omega$ the whole phase space $\mathbb{S}$ of the
Markov chain $\xi(t)$ constitutes the (only) communication class of this chain.
\end{lemma}

\proof Fix an environment $\omega$ and consider matrices
\[
\tilde{P}_{n}\overset{\mathrm{def}}{=}(I-R_{n})^{-1}P_{n},\ \tilde{Q}%
_{n}\overset{\mathrm{def}}{=}(I-R_{n})^{-1}Q_{n}.
\]
Remark that $\tilde{P}_{n}(i,j)$ is the probability that the random walk
${\xi}$ starting at $(n,i)$ would reach $(n+1,j)$ at the time of its first
exit from layer $n$; the probabilistic meaning of $\tilde{Q}_{n}(i,j)$ is
defined similarly. $\tilde{P}_{n}(i,j)\ge\varepsilon>0$ and $\tilde{Q}%
_{n}(i,j)\ge\varepsilon>0$ because of condition $\mathbf{C2}$. It is now
obvious that a random walk $\xi(\cdot)$ starting from any $z\in\mathbb{S}$
would reach any $z_{1}\in\mathbb{S}$ with a positive probability. $\Box$

Matrices of the form $(I-R-Q\psi)^{-1}$, $(I-R-Q\psi)^{-1}P$, and
$(I-R-Q\psi)^{-1}Q$ arise in the proofs of many statements below. We shall
list several elementary properties of these matrices.

\begin{lemma}
\label{lemma2} If condition $\mathbf{C2}$ is satisfied, $(P,Q,R)\in
\mathcal{J}_{0}$ and $\psi$ is any stochastic matrix, then there is
a constant $C$ depending only on $\varepsilon$ and $m$ such that
\begin{equation}
\left\|  (I-R-Q\psi)^{-1}\right\|  \le C. \label{est1}%
\end{equation}

\end{lemma}

\proof Note first that $||R^l||\le 1-\varepsilon$ implies that for
some $C_1$ uniformly in $R$
\[
||(I-R)^{-1}||\le\sum_{k=0}^{\infty}||R^k||\le C_1.
\]
Next, it follows from $(P+Q+R)\mathbf{1}=\mathbf{1}$ that $(I-R)^{-1}%
P\mathbf{1}+(I-R)^{-1}Q\mathbf{1}=\mathbf{1}$ and $(I-R)^{-1}Q\mathbf{1}%
=\mathbf{1}-(I-R)^{-1}P\mathbf{1}$. Condition $\mathbf{C2}$ implies
that $(I-R)^{-1}P\mathbf{1}\ge m\varepsilon\mathbf{1}$. Hence
\[
\left\|  (I-R)^{-1}Q\right\|  =\left\|  (I-R)^{-1}Q\mathbf{1}\right\|
=\left\|  \mathbf{1}-(I-R)^{-1}P\mathbf{1}\right\|  \le1-m\varepsilon.
\]
Similarly, $\left\|  (I-R)^{-1}P\right\|  \le1-m\varepsilon$. Hence
\[
\begin{aligned}
\left\|(I-R-Q\psi)^{-1}\right\|&=\left\|
(I-(I-R)^{-1}Q\psi)^{-1}(I-R)^{-1}\right\|\\
&\le (1-\left\|(I-R)^{-1}Q\psi\right\|)^{-1}
\left\|(I-R)^{-1}\right\|\le C_1m^{-1}\varepsilon^{-1}\equiv C.
\end{aligned}
\]
Lemma is proved. $\Box$

\begin{lemma}
\label{lemma3} (\cite{BG}) If condition $\mathbf{C2}$ is satisfied,
$(P,Q,R)\in\mathcal{J}$, and $\psi$ is a stochastic matrix, then
$(I-R-Q\psi)^{-1}P$ is also stochastic.
\end{lemma}

\proof We have to check that $(I-R-Q\psi)^{-1}P\mathbf{1}=\mathbf{1}$ which is
equivalent to $P\mathbf{1}=(I-Q\psi-R)\mathbf{1}$ $\Leftrightarrow$
$(P+Q\psi+R)\mathbf{1}=\mathbf{1}.$ Since $\psi\mathbf{1}=\mathbf{1}$ and
$P+Q+R$ is stochastic, the result follows. $\Box$

\begin{lemma}
\label{lemma4} Suppose that condition $\mathbf{C2}$ is satisfied and
$(P,Q,R)\in\mathcal{J}_{0}$ and let a matrix $\varphi\ge0$ be such that
$\varphi\mathbf{1}\le\mathbf{1}$. Then
\begin{equation}
\label{est2}((I-R-Q\varphi)^{-1}P)(i,j)\ge\varepsilon
\ \hbox{ and }\ ((I-R-Q\varphi)^{-1}Q)(i,j)\ge\varepsilon.
\end{equation}

\end{lemma}

\proof $(I-R-Q\varphi)^{-1}P\ge(I-R)^{-1}P$ and $(I-R-Q\varphi)^{-1}%
Q\ge(I-R)^{-1}Q$. $\Box$

\subsection{ Random transformations, related Markov chains, Lyapunov
exponents, and recurrence criteria.\label{randomtransformations}}

The purpose of this section is to introduce objects listed in its title. These
objects shall play a major role in the proofs of our main results. They shall
also allow us to state the main results from \cite{BG} in the form which is
suitable for our purposes.

\smallskip\noindent\textbf{Random transformations and related Markov chains.}

\smallskip\noindent Let $\Psi$ be the set of stochastic $m\times m$ matrices,
$\mathbb{X}$ be the set of unit vectors with non-negative components, and
$\mathrm{M}\overset{\mathrm{def}}{=}\Psi\times\mathbb{X}$ the direct product
of these two sets. Define a distance $\rho(\cdot,\cdot)$ on $\mathrm{M}$ by
\begin{equation}
\rho((\psi,x),(\psi^{\prime},x^{\prime}))\overset{\mathrm{def}}{=} ||\psi
-\psi^{\prime}||+||x-x^{\prime}||. \label{distance}%
\end{equation}
For any triple $(P,Q,R)\in\mathcal{J}_{0}$ denote by $g\equiv g_{(P,Q,R)}$ a
transformation
\begin{equation}
\label{transf}g:\,\mathrm{M}\mapsto\mathrm{M,}\ \hbox{ where }\ g.(\psi
,x)\overset{\mathrm{def}}{=}((I-R-Q\psi)^{-1}P\,,\, ||Bx||^{-1}Bx),
\end{equation}
and
\begin{equation}
B\equiv B_{(P,Q,R)}(\psi)\overset{\mathrm{def}}{=}(I-R-Q\psi)^{-1}Q.
\label{DefinB}%
\end{equation}
The fact that $g$ maps $\mathrm{M}$ into itself follows from Lemma
\ref{lemma3}.

\smallskip\noindent\textit{Remark.} Here and in the sequel the notation
$g.(\psi,x)$ is used instead of $g((\psi,x))$ and the dot is meant to replace
the brackets and to emphasize the fact that $g$ maps $(\psi,x)$ into another
pair from $\mathrm{M}$. In fact this notation is often used in the theory of
products of random matrices, e. g. $B.x\overset{\mathrm{def}}{=}
||Bx||^{-1}Bx$; we thus have extended this tradition to another component of
$g$.

If $\omega\in\Omega$ is an environment, $\omega=(\omega_{n})_{-\infty
<n<\infty}$, where $\omega_{n}\overset{\mathrm{def}}{=}(P_{n},Q_{n},R_{n})
\in\mathcal{J}_{0}$, then (\ref{transf}) allows us to define a sequence
$g_{n}\equiv g_{\omega_{n}}$ of random transformations of $\mathrm{M}$. Given
the sequence $g_{n}$, we define a Markov chain with a state space
$\mathcal{J}_{0}\times\mathrm{M}$. To this end consider an $a\in\mathbb{Z},$
and a $(\psi_{a}, x_{a})\in\mathrm{M}$ and put for $n\geq a$
\begin{equation}
\label{EqPsi}(\psi_{n+1},x_{n+1})\overset{\mathrm{def}}{=} g_{n}.(\psi
_{n},x_{n})\equiv((I-R_{n}-Q_{n}\psi_{n})^{-1}P_{n}\,,\,\|B_{n}x_{n}%
\|^{-1}B_{n}x_{n}),
\end{equation}
where we use a concise notation for matrices defined by (\ref{DefinB}):
\begin{equation}
B_{n}\overset{\mathrm{def}}{=} B_{\omega_{n}}(\psi_{n})\equiv B_{(P_{n}%
,Q_{n},R_{n})}(\psi_{n}). \label{defB}%
\end{equation}

\begin{theorem}
\label{ThZeta} Suppose that Condition $\mathbf{C}$ is satisfied. Then:

\noindent a) For $\mathbb{P}$-a.e. sequence $\omega$ the following limits
exist:
\begin{equation}
\zeta_{n}\overset{\mathrm{def}}{=}\lim_{a\rightarrow-\infty}\psi_{n}%
,\ \ y_{n}\overset{\mathrm{def}}{=}\lim_{a\rightarrow-\infty}x_{n}.
\label{DefZeta}%
\end{equation}
and $(\zeta_{n},y_{n})$ does not depend on the choice of the sequence
$(\psi_{a},y_{a}).$ Furthermore, the convergence in (\ref{DefZeta}) is uniform
in $(\psi_{a},x_{a})$.

\noindent b) The sequence of pairs $(\zeta_{n},y_{n})\equiv(\zeta_{n}%
(\omega),y_{n}(\omega))\ -\infty<n<\infty,$ is the unique sequence of elements
from $\mathrm{M}$ which satisfy the following infinite system of equations
\begin{equation}
\label{EqZeta}(\zeta_{n+1},y_{n+1})=\left(  (I-R_{n}-Q_{n}\zeta_{n})^{-1}%
P_{n}\,,\, ||A_{n}(\omega)y_{n}||^{-1}A_{n}(\omega)y_{n}\right)  ,\quad
n\in\mathbb{Z},
\end{equation}
where
\begin{equation}
\label{DefA}A_{n}\equiv A_{n}(\omega)\overset{\mathrm{def}}{=} (I-R_{n}%
-Q_{n}\zeta_{n})^{-1}Q_{n}.
\end{equation}
c) The enlarged sequence $(\omega_{n},\zeta_{n},y_{n}),\ -\infty<n<\infty,$
forms a stationary and ergodic Markov chain with components $\omega_{n}$ and
$(\zeta_{n},y_{n})$ being independent of each other.
\end{theorem}

\proof  The first relation in (\ref{DefZeta}) is the most important
statement of our Theorem and it also is the main content of Theorem
1 in \cite{BG}; it thus is known.

The main difference between this Theorem and Theorem 1 from \cite{BG} is that
here we consider the extended sequence $(\psi_{n},x_{n}),\ n\ge a$, rather
than just $(\psi_{n}),\ n\ge a$. The proof of the second relation in
(\ref{DefZeta}) is based on two observations. First note that the first
relation in (\ref{DefZeta}) implies that $\lim_{a\rightarrow-\infty}%
B_{n}=A_{n}$. Next, it follows from the definition of the sequence $x_{n}$
that
\begin{equation}
\label{DefY}x_{n}=\left\|  B_{n-1}\dots B_{a}x_{a}\right\|  ^{-1} B_{n-1}\dots
B_{a}x_{a}.
\end{equation}
Estimates (\ref{est1}) and (\ref{est2}) imply that
$\min_{i_{1},i_{2}
,i_{3},i_{4}}B_{k}^{-1}(i_{1},i_{2})B_{k}(i_{3},i_{4})\ge\bar{\varepsilon}$
for some $\bar{\varepsilon}>0$ and hence also $
\min_{i_{1},i_{2},i_{3},i_{4}}A_{k}^{-1}(i_{1},i_{2})A_{k}(i_{3},i_{4})\ge
\bar{\varepsilon}. $ It is well known (and can be easily derived
from Lemma \ref{pos}) that these inequalities imply the existence of
\[
\lim_{a\rightarrow-\infty}\left\|  A_{n}A_{n-1}\dots A_{a}x_{a}\right\|  ^{-1}
A_{n}A_{n-1}\dots A_{a}x_{a}
\]
and this limit does not depend on the choice of the sequence $x_{a}\ge0,
||x_{a}||=1$. Combining these two limiting procedures we obtain the proof of
the second relation in (\ref{DefZeta}).

Part b) of the Theorem is proved exactly as part b) of Theorem 1 from
\cite{BG}.

The Markov chain property and the independence claimed in part c) are obvious
corollaries of the independence of the triples $(P_{n},Q_{n},R_{n})$. And,
finally, the ergodicity of the sequence $(\omega_{n}, \zeta_{n}, y_{n})$ is
due to the fact that the sequence $\omega_{n}$ is ergodic and the $(\zeta_{n},
y_{n})$ is a function of $(\omega_{k})_{k\le n-1}$. $\Box$

\textit{Remark.} The proof of Theorem 1 in \cite{BG} was obtained under much
less restrictive assumptions than those listed in condition $\mathbf{C}$ of
this work. In particular, the i. i. d. condition which we impose on our
environments (rather than having them just stationary and ergodic) is
unimportant for parts a) and b) of Theorem \ref{ThZeta} as well as for Theorem
\ref{ThLyap}. However, the i. i. d. property is important for the proof of our
main results.

\smallskip\noindent\textbf{The top Lyapunov exponent of products of matrices
$A_{n}$ and the recurrence criteria.}

\smallskip\noindent The top Lyapunov exponent of products of matrices $A_{n}$
will be denoted by $\lambda$ and it is defined by
\begin{equation}
\label{DefLyap}\lambda\overset{\mathrm{def}}{=}\lim_{n\rightarrow\infty}%
{\frac{1}{{n}}}\log\left\|  A_{n}A_{n-1}\dots A_{1}\right\|  .
\end{equation}
The existence of the limit in (\ref{DefLyap}) with $\mathbb{P}$-probability 1
and the fact that $\lambda$ does not depend on $\omega$ is an immediate
corollary of the Kingman's sub-additive ergodic theorem; it was first proved
in \cite{FK}. The Furstenberg formula states that
\begin{equation}
\label{LyapF}\lambda=\int_{\mathcal{J}_{0}\times\mathrm{M}}\log\left\|
(I-R-Q\zeta)^{-1}Qy\right\|  \mu(dg)\nu(d(\zeta,y)),
\end{equation}
where $\nu(d(\zeta,y))$ is the invariant measure of the Markov chain
(\ref{EqPsi}) and $\mu(dg)$ is the distribution of the set of
triples $(P,Q,R)$ supported by $\mathcal{J}_{0}$ (defined in section
\ref{mainresults}). We use the shorter notation $dg$ rather than
$d(P,Q,R)$ because, as we have seen above, every triple
$(P,Q,R)\in\mathcal{J}_{0}$ defines a transformation $g$. Besides,
this notation is consistent with the one used in section
\ref{rand-tr}.

We remark that a proof of (\ref{DefLyap}) and (\ref{LyapF}) will be given in
section \ref{sec2.4} as a natural part of the proof of the invariance
principle for the sequence of random variables $\log\left\|  A_{n}A_{n-1}\dots
A_{1}\right\|  $.

We finish this section by quoting the recurrence criteria proved in \cite{BG}.

\begin{theorem}
\label{ThLyap} Suppose that Condition $\mathbf{C}$ is satisfied. Then

\noindent a) $\lambda\gtrless0$ if and only if for $\mathbb{P}$-a.e.
environment $\omega$ one has (respectively)
\[
\lim_{t\rightarrow\infty}\xi(t)=\mp\infty
\ \ \hbox{$Pr_{\omega}$-almost surely}.
\]
b) $\lambda=0$ if and only if for $\mathbb{P}$-a.e. $\omega$ the RW $\xi
(\cdot)$ is recurrent, that is
\[
\limsup_{t\rightarrow\infty}\xi(t)=+\infty\ \hbox{ and }\ \liminf
_{t\rightarrow\infty} \xi(t)=-\infty\ \ \hbox{$Pr_{\omega}$-almost surely}.
\]

\end{theorem}

\subsection{One algebraic corollary of Theorems \ref{ThZeta} and
\ref{ThLyap}.\label{sec2.3}}

Theorems \ref{ThZeta} and \ref{ThLyap} combined with a simple probabilistic
observation lead to an algebraic result which plays a very important role in
the proof of our algebraic condition.

Suppose that the matrices $(P_{n},Q_{n},R_{n})$ do not depend on $n$:
$(P_{n},Q_{n},R_{n})\equiv(P,Q,R)$, and the triple $(P,Q,R)$ satisfies
condition $\mathbf{C2}$. In this case relations (\ref{DefZeta}) mean that
$\zeta_{n}=\zeta$ and $y_{n}=y$, where $\zeta$ is a unique stochastic matrix
and $y\ge0$ a unique unit vector such that
\begin{equation}
\label{2.14}\zeta=(I-R-Q\zeta)^{-1}P,\ \hbox{ and }\ Ay=e^{\lambda}y,
\end{equation}
where the matrix $A$ is defined by
\[
A\overset{\mathrm{def}}{=}(I-R-Q\zeta)^{-1}Q.
\]
Theorem \ref{ThLyap} now states that a random walk in a constant environment
is recurrent if $\lambda=0$, transient to the right if $\lambda<0$, and
transient to the left if $\lambda>0$.

But the fact that the random environment does not depend on $n$ allows one to
analyse the recurrence and transience properties of the random walk in a way
which is much more straightforward than the one offered by Theorems
\ref{ThZeta} and \ref{ThLyap}.

Namely, suppose that $\xi(t)=(X(t),Y(t))=(k,i)$. Then the conditional
probability $Pr\{\,Y(t)=j\,|\,\xi(t-1)=(k,i)\}=P(i,j)+Q(i,j)+R(i,j)$ does not
depend on $X(t-1)$ and thus the second coordinate of this walk is a Markov
chain with a state space $(1,...,m)$ and a transition matrix $P+Q+R$. Hence,
if $\pi=(\pi_{1},...\pi_{m})$ is a probability vector such that $\pi
(P+Q+R)=\pi$ then $\pi_{i}$ is the frequency of visits by the RW to the sites
$(\cdot,i)$ of the strip.

Consider next the displacement $\eta(t)\overset{\mathrm{def}}{=}
X(t)-X(t-1)$ of the coordinate $X$ of the walk which occurs between
times $t-1$ and $t$. The random variable $\eta(t)$ takes values 1,
-1, or 0 and the following conditional distribution of the pair
$(\eta(t),Y(t))$ is given by
$Pr\{\,(\eta(t),Y(t))=(1,j)\,|\,\xi(t-1)=(k,i)\}=P(i,j)$,
$Pr\{\,(\eta (t),Y(t))=(-1,j)\,|\,\xi(t-1)=(k,i)\}=Q(i,j)$, and
$Pr\{\,(\eta (t),Y(t))=(0,j)\,|\,\xi(t-1)=(k,i)\}=R(i,j)$. It is
essential that this distribution depends only on $i$ (and not on
$k$) and thus this pair forms a time-stationary Markov chain. Let us
denote by $E_{(k,i)}$ the corresponding conditional expectation with
conditioning on $(\eta(t-1), Y(t-1))=(k,i)$, $-1\le k\le1,\ 1\le m$.
We then have
\[
E_{(k,i)}(\eta(t))=\sum_{j=1}^{m}P(i,j)-\sum_{j=1}^{m}Q(i,j).
\]
and the expectation of the same random variable with respect to the stationary
distribution is thus given by $\sum_{i=1}^{m}\pi_{i}\sum_{j=1}^{m}%
(P(i,j)-Q(i,j))$. Applying the law of large numbers for Markov chains to the
sequence $\eta(t)$ we obtain that with $Pr$-probability 1
\[
\lim_{t\to\infty}t^{-1}X(t)=\lim_{t\to\infty}t^{-1}\sum_{k=1}^{t}\eta(k)
=\sum_{i=1}^{m}\pi_{i}\sum_{j=1}^{m}(P(i,j)-Q(i,j))
\]
and this limit is independent of the $\xi(0)$. Since this result is equivalent
to the statements of Theorems \ref{ThZeta} and \ref{ThLyap}, we obtain the following

\begin{lemma}
. \label{alg-equiv} Suppose that $(P,Q,R)$ satisfies Condition
$\mathbf{C2}$. Then $(\zeta,x)\in\mathrm{M}$ satisfies equations
(\ref{2.14}) with $\lambda=0$ if and only if
\begin{equation}
\sum_{i=1}^{m}\pi_{i}\sum_{j=1}^{m}(P(i,j)-Q(i,j))=0.
\end{equation}
Moreover $\lambda>0$ if and only if $\sum_{i=1}^{m}\pi_{i}\sum_{j=1}%
^{m}(P(i,j)-Q(i,j))<0$ (and thus $\lambda<0$ if and only if $\sum_{i=1}^{m}%
\pi_{i}\sum_{j=1}^{m}(P(i,j)-Q(i,j))>0$).
\end{lemma}

\subsection{ The CLT and the invariance principle for $S_{n}$'s.\label{sec2.4}%
}

The main goal of this section is to prove an invariance principle (IP) (and a
CLT) for the sequence
\begin{equation}
\label{Sn}S_{n}\overset{\mathrm{def}}{=}\log\left\|  B_{n}\dots B_{1}%
x_{1}\right\|  -n\lambda,
\end{equation}
where matrices $B_{n}$ are defined by (\ref{defB}) and $\lambda$ is given by
(\ref{LyapF}). Obviously, $S_{n}$ depends on $(\psi_{1},x_{1})\in\mathrm{M}$.
We shall prove that in fact the IP (and the CLT) are satisfied uniformly in
$(\psi_{1},x_{1})\in\mathrm{M}$. Moreover, exactly one of the two things takes
place if the random walk is recurrent: either the asymptotic behaviour of
$S_{n}$ is described by a non-degenerate Wiener process, or the support of the
distribution of matrices $(P,Q,R)$ belongs to an algebraic manifold defined by
equations (\ref{algebra}).

To make these statements precise we first recall one of the definitions of the
invariance principle associated with a general random sequence $S_{n}%
=\sum_{k=1}^{n}f_{k}$, with the convention $S_{0}=0$. Let $\{
C[0,1],\mathcal{B},P_{W}\}$ be the probability space where $C[0,1]$ is the
space of continuous functions with the $\mathrm{sup}$ norm topology,
$\mathcal{B}$ being the Borel $\sigma$-algebra generated by open sets in
$C[0,1]$, and $P_{W}$ the Wiener measure. Define for $t\in[0,1]$ a sequence of
random functions $v_{n}(t)$ associated with the sequence $S_{n}$. Namely, put
\begin{equation}
\label{defv}v_{n}(t)\overset{\mathrm{def}}{=} n^{-\frac{1}{2}}\left(  {S}%
_{k}+f_{k+1}(tn-k)\right)  \hbox{
if $k\le tn\le k+1,\ k=0,1,...,n-1$.}
\end{equation}
For a $\sigma>0$ let $\{\mathbb{P}_{n}^{\sigma}\}$ be the sequence of
probability measures on $\{\, C[0,1],\mathcal{B}\,\}$ determined by the
distribution of $\{\,\sigma^{-1}v_{n}(t),\ 0\le t\le1\,\}$.

\smallskip\noindent\textbf{Definition.} A random sequence $S_{n}$ satisfies
the invariance principle with parameter $\sigma>0$ if $\mathbb{P}_{n}^{\sigma
}\rightarrow P_{W}$ weakly as $n\to\infty$. If the sequence $S_{n}$ depends on
(another) parameter, e.g. $z_{1}$, then we say that $S_{n}$ satisfies the
invariance principle with parameter $\sigma>0$ uniformly in $z_{1}$ if for any
continuous functional on $\mathfrak{f}:C[0,1]\mapsto\mathbb{R}$ one has:
$\mathbb{E}_{n}^{\sigma}(\mathfrak{f})\rightarrow E_{W}(\mathfrak{f})$
uniformly in $z_{1}$ as $n\to\infty$. Here $\mathbb{E}_{n}$ and $E_{W}$ are
expectations with respect to the relevant probabilities.

Let us state the invariance principle for the sequence $S_{n}$ given by
(\ref{Sn}). Note that in this case
\begin{equation}
S_{n}=\sum_{k=1}^{n}(\log\left\|  B_{k}x_{k}\right\|  -\lambda), \hbox{
where } x_{k}=\left\|  B_{k-1}x_{k-1}\right\|  ^{-1}B_{k-1}x_{k-1},\ k\ge2.
\end{equation}
Put $z_{n}=(\psi_{n},x_{n})$ and $f_{n}=f(g_{n},z_{n})$, where the function
$f$ is defined on the set of pairs $(g,z)\equiv((P,Q,R),(\psi,x))$ by
\begin{equation}
\label{deff}f(g,z)\overset{\mathrm{def}}{=}\log\left\|  (I-R-Q\psi
)^{-1}Qx\right\|  -\lambda.
\end{equation}
Obviously in these notations $S_{n}=\sum_{k=1}^{n}f_{k}$. Denote by
$\mathfrak{A}$ the Markov operator associated with the Markov chain
$z_{n+1}=g_{n}.z_{n}$ defined by (\ref{EqPsi}): if $F$ is a function defined
on the state space $\mathcal{J}_{0}\times\mathrm{M}$ of this chain then
\[
(\mathfrak{A}F)(g,z)\overset{\mathrm{def}}{=} \int_{\mathcal{J}_{0}%
\times\mathrm{M}} F(g^{\prime},g.z)\mu(dg^{\prime}).
\]
Using these notations we write $\nu(dz)$ (rather than $\nu(d(\psi,x))$) for
the invariant measure of the chain $z_{n}$ and we denote by $\mathrm{M}%
_{0}\subset\mathrm{M}$ the support of $\nu(dz)$.

\begin{theorem}
\label{IPmain} Suppose that condition $\mathbf{C}$ is satisfied and the
function $f$ is defined by (\ref{deff}). Then:

\noindent(i) The equation
\begin{equation}
\label{markov11}F(g,z)-(\mathfrak{A}F)(g,z)= f(g,z)
\end{equation}
has a unique solution $F(g,z)$ which is continuous on $\mathcal{J}_{0}%
\times\mathrm{M_{0}}$ and
\[
\int_{\mathcal{J}_{0}\times\mathrm{M}} F(g,z)\mu(dg)\nu(dz) =0.
\]

Denote by
\[
\sigma^{2} =\int_{\mathcal{J}_{0}\times\mathrm{M}_{0}} (\mathfrak{A}%
F^{2}-(\mathfrak{A}F)^{2})(g,y)\mu(dg)\nu(dy)
\]
(ii) If $\sigma>0$ then $\frac{{S}_{n}}{\sigma\sqrt n}$ converges in law
towards the standard Gaussian distribution $N(0,1)$ and the sequence ${S}_{n}$
satisfies the invariance principle with parameter $\sigma$ uniformly in
$(\psi_{1},x_{1})\in\mathrm{M}$.

\noindent(iii) If $\sigma=0$, then the function $F(g,y)$ depends only on $y$
and for every $(g,y)\in\mathcal{J}_{0}\times\mathrm{M}_{0}$ one has
\begin{equation}
\label{markov6}f(g,y) = F(y) -F(g.y).
\end{equation}

\noindent(iv) If $\sigma=0$ and $\lambda=0$ then
\begin{equation}
\label{markov7}\mathcal{J}_{0}\subset\mathcal{J}_{al},
\end{equation}
with $\mathcal{J}_{al}$ given by (\ref{algebra}).
\end{theorem}

\proof Statements (i), (ii), and (iii) of our Theorem follow from Theorem
\ref{IPMC}. In order to be able to apply Theorem \ref{IPMC} we have to show
that the sequence of random transformations $g_{n}$ has the so called
contraction property. Lemma \ref{contraction} establishes this property.
Relation (\ref{markov7}) is then derived from (\ref{markov6}) and one more
general property of Markov chains generated by products of contracting
transformations (Lemma \ref{support}).

\begin{lemma}
\label{contraction} Suppose that condition $\mathbf{C}$ is satisfied and let
\[
(\psi_{n+1},x_{n+1})=g_{n}.(\psi_{n},x_{n}),\ \ (\psi_{n+1}^{\prime}%
,x_{n+1}^{\prime})=g_{n}.(\psi_{n}^{\prime},x_{n}^{\prime}), \ \ n\ge1,
\]
be two sequences from $\mathrm{M}$. Then there is a $c,\ 0\le c<1,$ such that
for any $(\psi_{1},x_{1})$, $(\psi_{1}^{\prime},x_{1}^{\prime})\in\mathrm{M}$
\begin{equation}
\rho\left(  (\psi_{n},x_{n}),(\psi_{n}^{\prime},x_{n}^{\prime})\right)
\le\mathrm{const}\,c^{n},
\end{equation}
where $\rho(\cdot,\cdot)$ is defined by (\ref{distance}).
\end{lemma}

\textbf{Proof} of Lemma \ref{contraction}. We shall first prove that there is
a $c_{0}<1$ such that $||\psi_{n}-\psi_{n}^{\prime}||\le\mathrm{const}%
\,c_{0}^{n}$. The control of the $x$-component would then follow from this result.

Let us introduce a sequence of $m\times m$ matrices $\varphi_{n}$, $n\ge1$,
which we define recursively: $\varphi_{1}=0$ and
\begin{equation}
\label{phi1}\varphi_{n+1}=(I-R_{n}-Q_{n}\varphi_{n})^{-1}P_{n}%
,\ \hbox{ if $n\ge 1$}.
\end{equation}
\textit{Remark.} Matrices $\varphi_{n}$ and $\psi_{n}$ were defined in a
purely analytic way. Their probabilistic meaning is well known (see \cite{BG})
and shall also be discussed in Section \ref{sec3.1}.

Put $\Delta_{k}\overset{\mathrm{def}}{=}\psi_{k}-\varphi_{k}$. To control the
$\psi$-part of the sequence $(\psi_{n},x_{n})$ we need the following

\begin{lemma}
\label{contrpsi} Suppose that condition $\mathbf{C}$ is satisfied. Then there
is a $c_{0},\ 0\le c_{0}<1,$ such that for any stochastic matrix $\psi_{1}%
\in\Psi$ the matrix elements of the corresponding $\Delta_{n+1}$ are of the
following form:
\begin{equation}
\label{EqD0}{\Delta_{n+1}(i,j)}=\alpha_{n}(i)c_{n}(j)+\tilde{\epsilon}%
_{n}(i,j).
\end{equation}
Here $\alpha_{n}(i)$ and $c_{n}(j)$ depend only on the sequence $(P_{j}%
,Q_{j},R_{j}), 1\le j\le n$;

\noindent the matrix $\tilde{\epsilon}_{n}=(\tilde{\epsilon}_{n}(i,j))$ is a
function of $\psi_{1}$ and of the sequence $(P_{j},Q_{j},R_{j})$, $1 \le j \le
n,$ satisfying $||\tilde{\epsilon}_{n}||\le C_{1}c_{0}^{n}$ for some constant
$C_{1}$.
\end{lemma}

\textbf{Corollary.} If Condition $\mathbf{C}$ holds then
\begin{equation}
\label{controlpsi}||\psi_{n+1}-\psi_{n+1}^{\prime}||\le2C_{1}\,c_{0}^{n}.
\end{equation}
\textbf{Proof} of Corollary. Consider a sequence $\psi_{n}^{\prime}$ which
differs from $\psi_{n}$ in that the starting value for recursion (\ref{EqPsi})
is $\psi_{1}^{\prime}$. Put $\Delta_{k}^{\prime}\overset{\mathrm{def}}{=}%
\psi_{k}^{\prime}-\varphi_{k}$. Applying the result of Lemma \ref{contrpsi} to
$\Delta_{n+1}^{\prime}$ we obtain:
\begin{equation}
\label{EqD4}{\Delta_{n+1}^{\prime}(i,j)}=\alpha_{n}(i)c_{n}(j)+\tilde
{\epsilon}^{\prime}_{n}(i,j).
\end{equation}
It follows from (\ref{EqD0}), (\ref{EqD4}), and the definition of
$\Delta_{n+1}$ and $\Delta_{n+1}^{\prime}$ that $||\psi_{n+1}-\psi
_{n+1}^{\prime}||=||\Delta_{n+1}-\Delta_{n+1}^{\prime}|| \le||\tilde{\epsilon
}_{n}||+||\tilde{\epsilon}^{\prime}_{n}||\le2C_{1}\,c_{0}^{n}$. $\Box$

\smallskip\noindent\textbf{Proof} of Lemma \ref{contrpsi}. The main idea of
this proof is the same as that of the proof of Theorem 1 from \cite{BG}. A
very minor difference is that here we have to control the behaviour of
$\psi_{n}$ when $n$ is growing while $\psi_{1}$ is fixed; in \cite{BG} $n$ was
fixed while the starting point of the chain was tending to $-\infty$. A more
important difference is that here we state the exponential speed of
convergence of certain sequences and present the corresponding quantities in a
relatively explicit way while in \cite{BG} the speed of convergence was not
very essential (even though the exponential character of convergence had been
clear already then).

To start, note that it follows from (\ref{EqPsi}) and (\ref{phi1}) that
\begin{equation}
\label{EqD}\begin{aligned}
\Delta_{n+1}&=((I-R_{n}-Q_{n}\psi_{n})^{-1}-(I-R_n-Q_n
\varphi_n)^{-1})P_{n}\\
&=(I-R_{n}-Q_{n}\psi_{n})^{-1}Q_n\Delta_n(I-R_n-Q_n\varphi_n)^{-1}P_{n}=
B_n\Delta_n\varphi_{n+1} \end{aligned}
\end{equation}
Iterating (\ref{EqD}), we obtain
\begin{equation}
\label{EqD0.1}\Delta_{n+1}= B_{n}...B_{1}\Delta_{1}\varphi_{2}...\varphi
_{n+1}\equiv B_{n}...B_{1}\psi_{1}\varphi_{2}...\varphi_{n+1}.
\end{equation}
It follows from Lemma \ref{lemma4} that
$\varphi_{n}\mathbf{1}\le\mathbf{1}$. The matrix elements of the
matrices $\varphi_{n}$, $n\ge2$, are strictly positive and,
moreover, according to estimates (\ref{est2}) we have:
$\varphi_{n}(i,j)\ge\varepsilon$ (and hence also
$\varphi_{n}(i,j)\le 1 -(m-1)\varepsilon$). We are in a position
to apply to the product of matrices $\varphi_{n}$ the presentation
derived in Lemma \ref{pos} (with $a_n$'s replaced by
$\varphi_{n}$'s). By the first formula in (\ref{A1}), we have:
\[
\varphi_{2}...\varphi_{n+1}=D_{n}[\left(  c_{n}(1)\mathbf{1},\ldots
,c_{n}(m)\mathbf{1}\right)  +\phi_{n}],
\]
where $D_{n}$ is a diagonal matrix, $c_{n}(j)\ge\delta$ with $\sum
_{j=1}^{m}c_{n}(j)=1$, and $\left\|  \phi_{n}\right\|
\le(1-m\delta)^{n-1}$ with $\delta>0$ (and of course $m\delta<1$).
One can easily see that $\delta\ge m^{-1}\varepsilon^2$ (this
follows from (\ref{delta}) and the above estimates for
$\varphi_{n}(i,j)$). We note also that the estimate for $c_{n}(j)$
follows from (\ref{delta1}) and (\ref{delta2}).

Put $c_{0}=1-m\delta$ and let
$\mathcal{B}_{n}\overset{\mathrm{def}}{=}
B_{n}...B_{1}\Delta_{1}D_{n}$. We then have
\begin{equation}
\label{EqD1}\Delta_{n+1}= \mathcal{B}_{n}[\left(  c_{n}(1)\mathbf{1}%
,\ldots,c_{n}(m)\mathbf{1}\right)  +\phi_{n}]
\end{equation}
and thus $\Delta_{n+1}(i,j)=c_{n}(j)\sum_{k=1}^{m}\mathcal{B}_{n}(i,k)(1
+\frac{\phi_{n}(k,j)}{c_{n}(j)})$. But all $\mathcal{B}_{n}(i,k)>0$ and
$\max_{k,j}|\phi_{n}(k,j)|{c^{-1}_{n}(j)} \le\mathrm{const}\, c_{0}^{n} $.
Hence
\begin{equation}
\label{EqD2}\frac{\Delta_{n+1}(i,l)}{\Delta_{n+1}(i,j)}= \frac{c_{n}(l)}%
{c_{n}(j)}+\epsilon_{n}(i,j,l),
\end{equation}
where $|\epsilon_{n}(i,j,l)|<Cc_{0}^{n}$ with $C$ being some constant. It
follows from (\ref{EqD2}) that
\[
(\Delta_{n+1}(i,j))^{-1}\sum_{l=1}^{m}\Delta_{n+1}(i,l)= \frac{1}{c_{n}%
(j)}+\epsilon_{n}(i,j).
\]
On the other hand remember that
\[
\sum_{l=1}^{m}\Delta_{n+1}(i,l)=\sum_{l=1}^{m}\psi_{n+1}(i,l)- \sum_{l=1}%
^{m}\varphi_{n+1}(i,l)= 1-\sum_{l=1}^{m}\varphi_{n+1}(i,l)\overset
{\mathrm{def}}{=} \alpha_{n}(i).
\]
Comparing these two expressions we obtain that
\begin{equation}
\label{EqD3}{\Delta_{n+1}(i,j)}=\alpha_{n}(i)c_{n}(j)+\tilde{\epsilon}%
_{n}(i,j),
\end{equation}
where $|\tilde{\epsilon}_{n}(i,j)|\le C_{1}c_{0}^{n}$. Lemma \ref{contrpsi} is
proved. $\Box$

We now turn to the difference $|| x_{n+1}-x_{n+1}^{\prime}||$. Let us denote
by $b_{n}$ the transformation of the set $\mathbb{X}$ of unit non-negative
vectors defined by
\begin{equation}
\label{defb}b_{n}(x)=||B_{n}x||^{-1}B_{n}x,\ \hbox{ where
$B_n=(I-R_n-Q_n\psi_n)^{-1}Q_n$},
\end{equation}
and $\psi_{n}$ are the same as above. The sequence $b_{n}^{\prime}$ is defined
in a similar way with the only difference that $\psi_{n}$ is replaced by
$\psi_{n}^{\prime}$. Inequality (\ref{controlpsi}) implies that for some
$C_{2}$
\[
\bar{\rho}(b_{n},b_{n}^{\prime})\overset{\mathrm{def}}{=} \sup_{x\in
\mathbb{X}}||b_{n}(x)-b_{n}^{\prime}(x)||\le C_{2}c_{0}^{n}.
\]
A very general and simple Lemma \ref{stability} from Appendix now implies
that
\[
|| x_{n+1}-x_{n+1}^{\prime}||\le C(\epsilon)(c_{0}+\epsilon)^{n}(1+||
x_{1}-x_{1}^{\prime}||)
\]
and this proves Lemma \ref{contraction}. $\Box$

We can now easily prove the existence of the limit in (\ref{DefLyap}) as well
as the Furstenberg's formula (\ref{LyapF}) for $\lambda$. To this end note
that
\begin{equation}
\label{defbarS}\bar{S}_{n}(\zeta_{1},\mathbf{1})\overset{\mathrm{def}}{=}
\log||A_{n}...A_{1}||=\log||A_{n}...A_{1}\mathbf{1}||= \sum_{k=1}^{n}%
f(g_{k},z_{k})
\end{equation}
where the notation is chosen so that to emphasize the dependence of
the sum $\bar{S}_{n}(\zeta_{1},\mathbf{1})$ on initial values
$x_{1}=\mathbf{1}$ and $\psi_{1}=\zeta_{1}$ of the Markov chain.
(Remark the difference between $\bar{S}_{n}(\zeta_{1},\mathbf{1})$
and the sum $S_{n}$ in (\ref{Sn}).) Lemma \ref{contraction} implies
that
\begin{equation}
\label{differenceS}|\bar{S}_{n}(\zeta_{1},\mathbf{1}) -\bar{S}_{n}(\psi
_{1},x_{1})|\le C_{3},
\end{equation}
where the constant $C_{3}$ depends only on the parameter $\varepsilon$ from
condition $\mathbf{C}$. But then, according to the law of large numbers
applied to the Markov chain $(\omega_{n},\zeta_{n},y_{n})\equiv(g_{n}%
,\zeta_{n},y_{n}) $ defined in Theorem \ref{ThZeta} we have that the following
limit exists with probability 1:
\[
\lim_{n\to\infty}\frac{1}{n}\log||A_{n}...A_{1}||= \lim_{n\to\infty}\frac
{1}{n}\bar{S}_{n}(\zeta_{1},y_{1}) =\lambda,
\]
where $\lambda$ is given by (\ref{LyapF}).

Formula (\ref{LyapF}) implies that the mean value of the function
$f(g,z)$ defined by (\ref{deff}) is $0$. Also, it is obvious that
this function is Lipschitz on $\mathcal{J}_{0}\times\mathrm{M}$ in
all variables. Hence, Theorem \ref{IPMC} applies to the sequence
$S_{n}$ and statements (i), (ii), and (iii) of Theorem
\ref{IPmain} are thus proved.

\smallskip\noindent\textbf{The case $\sigma=0$ and $\lambda=0$: derivation of
the algebraic condition for $(P,Q,R)$.} We start with a statement which is a
corollary of a very general property proved in Lemma \ref{support0} from Appendix.

\begin{lemma}
\label{support} Suppose that Condition $\mathbf{C}$ is satisfied and let
$g\in\mathcal{J}_{0}$, $z_{g}\in\mathrm{M}$ be such that $g.z_{g}=z_{g}$. Then
$z_{g}\in\mathrm{M}_{0}\equiv\mathrm{supp}\nu$.
\end{lemma}

\proof According to Lemma \ref{contraction}, Condition $\mathbf{C}$ implies
that every $g\in\mathcal{J}_{0}$ is contracting. Hence, by Lemma
\ref{support0}, $z_{g}\in\mathrm{M}_{0}$. $\Box$

\textit{Derivation of the algebraic condition.} According to Theorem
\ref{IPMC} (see formula (\ref{markov5})), the equality $\sigma=0$
implies that $f(g,z) = F(z) -F(g.z)$. Hence, if $z$ can be chosen to
be equal to $z_{g}$, then it follows that $f(g,z_{g})=0$.

In the context of the present Theorem the function $f$ is given by $f(g,z) =
\log||(I-R-Q\psi)^{-1}Qx||$, where $g=(P,Q,R)\in\mathcal{J}_{0}$ and
$z=(\psi,x)\in\mathrm{M}_{0}\subset\Psi\times\mathbb{X}$. The equation
$g.z_{g}=z_{g}$ is equivalent to saying that $z_{g}=(\psi,x)$ satisfies
\[
(I-R-Q\psi)^{-1}\psi=\psi\ \ \hbox{and}\ \ ||(I-R-Q\psi)^{-1}Qx||^{-1}%
(I-R-Q\psi)^{-1}Qx=x.
\]
The equation $f(g,z_{g})=0$ now reads $\log||(I-R-Q\psi)^{-1}Qx||=0$ or,
equivalently, $||(I-R-Q\psi)^{-1}Qx||=1$. Hence the conditions $\sigma=0$ and
$\lambda=0$ imply that all pairs $(g,z_{g})\in\mathcal{J}_{0}\times
\mathrm{M}_{0}$ satisfy
\[
(I-R-Q\psi)^{-1}P=\psi\ \ \hbox{and}\ \ (I-R-Q\psi)^{-1}Qx=x.
\]
But, by Lemma \ref{alg-equiv}, this implies that $\mathcal{J}%
_{0}\subset\mathcal{J}_{al}$, where $\mathcal{J}_{al}$ is defined by
(\ref{algebra}). $\Box$

\section{Proof of Theorem \ref{Theor1}\label{sec3.1}}

As we are in the recurrent situation, we have that the Lyapunov exponent
$\lambda=0$.

Throughout this section we denote by $C$ a generic positive constant which
depends on nothing but $\varepsilon$ and $m$ and which may vary from place to
place. If $f,g>0$ are two functions, depending on $n\in\mathbb{Z}%
,\ i\in\left\{  1,\ldots,m\right\}  $, and maybe on other parameters, we write%
\[
f\asymp g\ \hbox{ if there exists a $C>1$ such that }\ C^{-1}f\leq g\leq Cf
\]

\textit{Potential and its properties.} As before, $S_{n}$ is defined by
(\ref{Sn}). We put
\begin{equation}
\Phi_{n}(\omega)\equiv\Phi_{n}\overset{\mathrm{def}}{=}\left\{
\begin{array}
[c]{ll}%
\log||A_{n}...A_{1}|| & \mathrm{if\ }n\geq1\\
0 & \mathrm{if\ }n=0\\
-\log||A_{0}...A_{n+1}|| & \mathrm{if\ }n\leq-1
\end{array}
\right.  \label{defPhi}%
\end{equation}
where the matrices $A_{n}$ are defined in (\ref{DefA}). If $n\geq1$, then
obviously $\Phi_{n}\equiv\bar{S}_{n}(\zeta_{1},\mathbf{1})$ defined in
(\ref{defbarS}). The random function $\Phi_{n}$ is the analog of the potential
considered first in \cite{S}. For $n\geq a$, $a\in\mathbb{Z},$ put
\begin{equation}
S_{a,n}(\omega;\psi_{a},x_{a})\equiv S_{a,n}(\omega)\overset{\mathrm{def}}%
{=}\log\left\Vert B_{n}\dots B_{a}x_{a}\right\Vert , \label{Sn1}%
\end{equation}
where the matrices $B_{n}$ are defined by (\ref{defB}). Similarly to
(\ref{differenceS}), one has that
\begin{equation}
\left\vert {S}_{a,n}(\omega;\zeta_{a},\mathbf{1})-{S}_{a,n}(\omega;\psi
_{a},x_{a})\right\vert \leq C, \label{differenceS1}%
\end{equation}
which implies:%
\begin{equation}
\left\vert S_{a,n}(\omega)-\left(  \Phi_{n}(\omega)-\Phi_{a}(\omega)\right)
\right\vert \leq C. \label{Sn2}%
\end{equation}

Since one of the conditions of Theorem \ref{Theor1} is $\mathcal{J}%
_{0}\not \subset \mathcal{J}_{al}$, it follows from Theorem \ref{IPmain}, part
\textit{(iv)} that $\Phi_{n}$ satisfies the invariance principle with a
strictly positive parameter $\sigma:\ \sigma>0$.

The importance of the potential $\left\{  \Phi_{n}\right\}
_{n\in\mathbb{Z}}$ is due to that fact that  it governs the
stationary measure of our Markov chain; in fact it defines this
stationary measure up to a multiplication by a bounded function (see
(\ref{Stat&Potential}). Namely, if $a<b,$ we consider the Markov
chain $\left\{ \xi_{t}^{a,b}\right\} _{t\in\mathbb{N}}$ on
\begin{equation}
\mathbb{S}_{a,b}\overset{\mathrm{def}}{=}\left\{  a,\ldots,b\right\}
\times\left\{  1,\ldots,m\right\}  \label{Def_Sab}%
\end{equation}
with transition probabilities (\ref{striptransition}) and reflecting
boundary conditions at $L_{a}$ and $L_{b}.$ This means that we
replace $\left( P_{a},Q_{a},R_{a}\right) $ by $\left(  I,0,0\right)
$ and $\left( P_{b},Q_{b},R_{b}\right) $ by $\left(  0,I,0\right) $.
This reflecting chain has a unique stationary probability measure
which we denote by $\pi_{a,b}=\left( \pi_{a,b}\left( k,i\right)
\right) _{\left( k,i\right) \in\mathbb{S}_{a,b}}.$ A description of
this measure was given in \cite{BG}. We repeat it here for the
convenience of the reader. To this end introduce row vectors
$\nu_k\de Z\left( \pi_{a,b}\left( k,i\right) \right) _{1\le i\le
m}$, $a\le k\le b$, and $Z$ is a (normalizing) factor. In terms of
these vectors the invariant measure equation reads
\begin{equation}
\begin{aligned}\label{invmeas}
& \nu_{k}=\nu_{k-1}P_{k-1}+\nu_{k}R_{k}+\nu_{k+1}Q_{k+1},\ \hbox{ if
}\ a<k<b \\
& \nu_a=\nu_{a+1}Q_{a+1},\ \ \nu_b=\nu_{b-1}P_{b-1}.
\end{aligned}
\end{equation}
To solve equations (\ref{invmeas}), define for $a\leq
k<b$ matrices $\alpha_{k}$ by%
\[
\alpha_{a}\overset{\mathrm{def}}{=}Q_{a+1},\ \ \hbox{and}\ \
\alpha_{k}\overset{\mathrm{def}}{=}Q_{k+1}\left(  I-R_{k}-Q_{k}\psi
_{k}\right)  ^{-1},\ \hbox{ when } \  a<k<b,
\]
where $\left\{  \psi_{k}\right\}  _{k\geq a+1}$ are given by
(\ref{EqPsi}) with the initial condition $\psi_{a+1}=I$ (we take
into account that $R_{a}=Q_{a}=0$ in our case). We shall now check
that $\nu_{k}$ can be found recursively as follows: $\nu_{k}
{=}\nu_{k+1}\alpha_{k},\ a\leq k<b,$, where $\nu_b$ satisfies
$\nu_b\psi_{b}=\nu_{b}$. Indeed, the boundary condition at $b$ in
(\ref{invmeas}) reduces to
$\nu_b=\nu_{b}\alpha_{b-1}P_{b-1}=\nu_b\psi_{b}$, where we use the
fact that $\alpha_{b-1}P_{b-1}=\psi_{b}$ because $Q_{b}=I$ (and also
due to (\ref{EqPsi})). But $\psi_{b}$ is an irreducible stochastic
matrix and therefore $\nu_{b}>0$ exists and is uniquely defined up
to a multiplication by a constant.  We now have for $a<k<b$ that
\begin{align*}
\nu_{k-1}P_{k-1}+\nu_{k}R_{k}+\nu_{k+1}Q_{k+1} &  =\nu_{k+1}\left(  \alpha
_{k}\alpha_{k-1}P_{k-1}+\alpha_{k}R_{k}+Q_{k+1}\right)  \\
&  =\nu_{k+1}\alpha_{k}\left(  Q_{k}\psi_{k}+R_{k}+\left(  I-R_{k}-Q_{k}%
\psi_{k}\right)  \right)  \\
&  =\nu_{k+1}\alpha_{k}=\nu_{k}.
\end{align*}
Finally $\nu_{a}=\nu_{a+1}Q_{a+1}$ with $\alpha_{a}=Q_{a+1}$ and
this finishes the proof of our statement.

We now have that
\[
\pi_{a,b}\left(  k,\cdot\right)  =\pi_{a,b}\left(  b,\cdot\right)
\alpha_{b-1}\alpha_{b-2}\cdot\cdots\cdot\alpha_{k},
\]
where as before $\pi_{a,b}\left(  k,\cdot\right)  $ is a row
vector. Note next that
\[
\alpha_{b-1}\alpha_{b-2}\cdot\cdots\cdot\alpha_{k}=B_{b-1}\cdot\cdots\cdot
B_{k+1}\left(  I-R_{k}-Q_{k}\psi_{k}\right)  ^{-1}.
\]
From this, we get%
\[
\pi_{a,b}\left(  k,\cdot\right)  \asymp\left\Vert B_{b-1}\cdot\cdots\cdot
B_{k+1}\right\Vert \pi_{a,b}\left(  b,\cdot\right)  ,
\]
and using (\ref{Sn1}), (\ref{Sn2}), we obtain for $a\leq k,l\leq b$%

\begin{equation}
\frac{\pi_{a,b}\left(  k,\cdot\right)  }{\pi_{a,b}\left(  l,\cdot\right)
}\asymp\exp\left[  \Phi_{k}-\Phi_{l}\right]  . \label{Stat&Potential}%
\end{equation}

We also consider the \textquotedblleft mirror
situation\textquotedblright\ by defining for $n\leq a$ the martices
$\psi_{n}^{-}$ in a similar way as in
(\ref{EqPsi}) by setting%
\[
\psi_{n-1}^{-}=\left(  I-R_{n}-P_{n}\psi_{n}^{-}\right)
^{-1}Q_{n},\ n\leq a,
\]
and a boundary condition $\psi_{a}^{-}$. Then, as in Theorem
\ref{ThZeta} a), one has that
$\zeta_{n}^{-}\overset{\mathrm{def}}{=}\lim_{a\rightarrow\infty
}\psi_{n}^{-}$ exists almost surely, and does not depend on the
boundary
condition $\psi_{a}^{-}$. We then put%
\[
A_{n}^{-}\overset{\mathrm{def}}{=}\left(
I-R_{n}-P_{n}\zeta_{n}^{-}\right) ^{-1}P_{n},
\]
and the potential $\Phi_{n}^{-}$ as (\ref{defPhi}):%
\[
\Phi_{n}^{-}\overset{\mathrm{def}}{=}\left\{
\begin{array}
[c]{ll}%
\log||A_{0}^{-}...A_{n-1}^{-}|| & \mathrm{if\ }n\geq1\\
0 & \mathrm{if\ }n=0\\
-\log||A_{n}^{-}...A_{-1}^{-}|| & \mathrm{if\ }n\leq-1
\end{array}
\right.  .
\]
We could as well have worked with this potential, and therefore\ we
obtain
\[
\frac{\pi_{a,b}\left(  k,\cdot\right)  }{\pi_{a,b}\left(  l,\cdot\right)
}\asymp\exp\left[  \Phi_{k}^{-}-\Phi_{l}^{-}\right]  .
\]
As $\Phi_{0}=\Phi_{0}^{-}=0,$ we get%
\begin{equation}
\left\vert \Phi_{n}-\Phi_{n}^{-}\right\vert \leq C\label{Compare_Phi_Phiprime}%
\end{equation}
uniformly in $n.$

It is convenient to slightly reformulate the invariance principle for the
potential. For that consider $C_{0}\left(  -\infty,\infty\right)  ,$ the space
of continuous functions $f:\left(  -\infty,\infty\right)  \rightarrow
\mathbb{R}$ satisfying $f\left(  0\right)  =0.$ We equip $C_{0}\left(
-\infty,\infty\right)  $ with a metric for uniform convergence on compacta,
e.g.%
\begin{equation}
d\left(  f,g\right)  \overset{\mathrm{def}}{=}\sum_{k=1}^{\infty}2^{-k}%
\min\left[  1,\sup\nolimits_{x\in\left[  -k,k\right]  }\left\vert f\left(
x\right)  -g\left(  x\right)  \right\vert \right]  ,\label{metric}%
\end{equation}
and write $\mathcal{B}$ for the Borel-$\sigma$-field which is also the
$\sigma$-field generated by the evaluation mappings $C_{0}\left(
-\infty,\infty\right)  \rightarrow\mathbb{R}.$ We also write $P_{W}$ for the
law of the double-sided Wiener measure on $C_{0}\left(  -\infty,\infty\right)
.$

For $n\in\mathbb{N},$ we define%
\[
W_{n}\left(  \frac{\left[  k\sigma^{2}\right]  }{n}\right)  \overset
{\mathrm{def}}{=}\frac{\Phi_{k}}{\sqrt{n}},\ k\in\mathbb{Z},
\]
and define $W_{n}\left(  t\right)  ,\ t\in\mathbb{R},$ by linear
interpolation. $W_{n}$ is a random variable taking values in $C_{0}\left(
-\infty,\infty\right)  .$

Weak convergence of $\left\{  W_{n}\left(  t\right)  \right\}  _{t\in
\mathbb{R}}$ on $C_{0}\left(  -\infty,\infty\right)  $ is the same as weak
convergence of $\left\{  W_{n}\left(  t\right)  \right\}  _{t\in\left[
-N,N\right]  }$ for any $N\in\mathbb{N},$ and therefore, we immediately get

\begin{proposition}
$W_{n}$ converges in law to $P_{W}.$
\end{proposition}

Let $V$ be the subset of functions $f\in C_{0}\left(
-\infty,\infty\right)  $ for which there exist real numbers
$a<b<c$ satisfying

\begin{enumerate}
\item
\[
0\in\left(  a,c\right)  .
\]

\item
\[
f\left(  a\right)  -f\left(  b\right)  =f\left(  c\right)  -f\left(  b\right)
=1
\]

\item
\[
f\left(  a\right)  >f\left(  x\right)  >f\left(  b\right)  ,\ \forall
x\in\left(  a,b\right)  ,
\]%
\[
f\left(  c\right)  >f\left(  x\right)  >f\left(  b\right)  ,\ \forall
x\in\left(  b,c\right)  .
\]

\item For any $\gamma>0$%
\begin{align*}
\sup_{x\in\left(  a-\gamma,a\right)  }f\left(  x\right)   &  >f\left(
a\right)  ,\\
\sup_{x\in\left(  c,c+\gamma\right)  }f\left(  x\right)   &  >f\left(
c\right)  .
\end{align*}

\end{enumerate}

It is clear that for $f\in V,$ $a,b,c$ are uniquely defined by $f,$ and we
write occasionally $a\left(  f\right)  ,b\left(  f\right)  ,c\left(  f\right)
$. $f\left(  b\right)  $ is the unique minimum of $f$ in $\left[  a,c\right]
.$ It is easy to prove that $V\in\mathcal{B},$ and%
\[
P_{W}\left(  V\right)  =1.
\]

If $\delta>0$ and $f\in V,$ we define%
\begin{align*}
&  c_{\delta}\left(  f\right)  \overset{\mathrm{def}}{=}\inf\left\{
x>c:f\left(  x\right)  =f\left(  c\right)  +\delta\right\} \\
&  a_{\delta}\left(  f\right)  \overset{\mathrm{def}}{=}\sup\left\{
x<a:f\left(  x\right)  =f\left(  a\right)  +\delta\right\}
\end{align*}
If $\gamma>0,$ we set $V_{\delta,\gamma}$ to be the set of functions $f\in V$
such that

\begin{enumerate}
\item
\begin{equation}
c_{\delta}\left(  f\right)  \leq1/\delta,\ a_{\delta}\left(  f\right)
\geq-1/\delta. \label{Est1_Delta}%
\end{equation}

\item
\begin{align}
\sup_{b\leq x<y\leq c_{\delta}}\left[  f\left(  x\right)  -f\left(  y\right)
\right]   &  \leq1-\delta,\label{Est2_Delta}\\
\sup_{a_{\delta}\leq y<x\leq b}\left[  f\left(  x\right)  -f\left(  y\right)
\right]   &  \leq1-\delta. \label{Est3_Delta}%
\end{align}

\item
\begin{equation}
\inf_{x\in\left[  a_{\delta},c_{\delta}\right]  \backslash\left(
b-\gamma,b+\gamma\right)  }f\left(  x\right)  \geq f\left(  b\right)  +\delta.
\label{Est4_Delta}%
\end{equation}

\end{enumerate}

It is evident that for any $\gamma>0,$ we have $V_{\delta,\gamma}\uparrow V$
for $\delta\downarrow0,$ and therefore, for any $\delta,\eta>0$ we can find
$\delta_{0}\left(  \gamma,\eta\right)  $ such that for $\delta\leq\delta_{0}$%
\[
P_{W}\left(  V_{\delta,\gamma}\right)  \geq1-\eta.
\]
It is easy to see that%
\[
P_{W}\left(  \partial V_{\delta,\gamma}\right)  =0,
\]
where $\partial$ refers to the boundary in $C_{0}\left(  -\infty
,\infty\right)  .$ Therefore, given $\gamma,\eta>0,$ we can find $N_{0}\left(
\gamma,\eta\right)  $ such that for $n\geq N_{0},$ $\delta\leq\delta_{0},$ we
have%
\begin{equation}
\mathbb{P}\left(  W_{n}\in V_{\delta,\gamma}\right)  \geq1-2\eta.
\label{Est_Yn}%
\end{equation}

For $t\in\mathbb{N},$ we set $n=n\left(  t\right)  \overset{\mathrm{def}}%
{=}\left[  \log^{2}t\right]  .$ If $W_{n\left(  t\right)  }\in V_{\delta
,\gamma}$ then we put%
\[
b_{t}\overset{\mathrm{def}}{=}\frac{b\left(  W_{n\left(  t\right)  }\right)
\log^{2}t}{\sigma^{2}},\ a_{t}\overset{\mathrm{def}}{=}\frac{a_{\delta}\left(
W_{n\left(  t\right)  }\right)  \log^{2}t}{\sigma^{2}},\ c_{t}\overset
{\mathrm{def}}{=}\frac{c_{\delta}\left(  W_{n\left(  t\right)  }\right)
\log^{2}t}{\sigma^{2}}.
\]

Remark that on $\left\{  W_{n\left(  t\right)  }\in V_{\delta,\gamma}\right\}
,$ we have the following properties, translated from (\ref{Est1_Delta}%
)-(\ref{Est4_Delta}):%
\begin{equation}
c_{t}\leq\frac{\log^{2}t}{\sigma^{2}\delta},\ a_{t}\geq-\frac{\log^{2}%
t}{\sigma^{2}\delta}, \label{Est1_Phi}%
\end{equation}%
\begin{align}
\Phi_{s}-\Phi_{s^{\prime}}  &  \leq\left(  1-\delta\right)  \log t,\ b_{t}\leq
s<s^{\prime}\leq c_{t},\label{Est2_Phi}\\
\Phi_{s}-\Phi_{s^{\prime}}  &  \leq\left(  1-\delta\right)  \log t,\ a_{t}\leq
s^{\prime}<s\leq b_{t}, \label{Est3_Phi}%
\end{align}%
\begin{equation}
\Phi_{s}\geq\Phi_{b_{t}}+\delta\log t,\ s\in\left[  a_{t},c_{t}\right]
\backslash\left[  b_{t}-\gamma\log^{2}t,b_{t}+\gamma\log^{2}t\right]  ,
\label{Est4_Phi}%
\end{equation}%
\begin{equation}
\min\left(  \Phi_{a_{t}},\Phi_{c_{t}}\right)  -\Phi_{b_{t}}\geq\left(
1+\delta\right)  \log t. \label{Est5_Phi}%
\end{equation}
Furthermore, if $0\in\left[  a_{t},b_{t}\right]  ,$ then%
\begin{equation}
\sup_{0\leq s\leq b_{t}}\Phi_{s}-\Phi_{b_{t}}\leq\log t, \label{Est6_Phi}%
\end{equation}
and similarly if $0\in\left[  b_{t},c_{t}\right]  .$

(We neglect the trivial issue that $a_{t},b_{t},c_{t}$ may not be in
$\mathbb{Z}$). The main result is

\begin{proposition}
\label{Prop_Main}For $\omega\in\left\{  W_{n\left(  t\right)  }\in
V_{\delta,\gamma}\right\}  ,$ we have for any $i\in\left\{  1,\ldots
,m\right\}  $%
\[
Pr_{\omega,\left(  0,i\right)  }\left(  X\left(  t\right)  \notin\left[
b_{t}-\gamma\log^{2}t,b_{t}+\gamma\log^{2}t\right]  \right)  \leq
4t^{-\delta/2},
\]
if $t$ is large enough.
\end{proposition}

Together with (\ref{Est_Yn}), this proves our main result Theorem \ref{Theor1}.

In all what follows, we keep $\gamma,\delta$ fixed, and assume that $\omega
\in\left\{  W_{n\left(  t\right)  }\in V_{\delta,\gamma}\right\}  .$ We will
also suppress $\omega$ in the notation, and will take $t$ large enough,
according to ensuing necessities.

We first prove several estimates of probabilities characterizing the behaviour
of a RW in a finite box in terms of the properties of the function $S_{n}$.

\begin{lemma}
\label{Le_Hittingtime1}Consider a random walk on $\mathbb{S}_{a,b}$ with
reflecting boundary conditions (see the discussion around (\ref{Def_Sab})),
and let $a<k<b.$ Then
\begin{align}
Pr_{\left(  k,i\right)  }\left(  \tau_{a}<\tau_{b}\right)   &  \leq
C\sum_{y=k}^{b}\exp\left(  \Phi_{y}-\Phi_{a}\right)  ,\label{hit1_1}\\
Pr_{\left(  k,i\right)  }\left(  \tau_{b}<\tau_{a}\right)   &  \leq
C\sum_{y=a}^{k}\exp\left(  \Phi_{y}-\Phi_{a}\right)  . \label{hit1_2}%
\end{align}
Here $\tau_{a},\tau_{b}$ are the hitting times of the layers $L_{a},L_{b}.$
\end{lemma}

\noindent\proof We only have to prove (\ref{hit1_1}). (\ref{hit1_2}) then
follows in the mirrored situation and using (\ref{Compare_Phi_Phiprime}).

Put $h_{k}(i)=Pr_{(k,i)}\left(  \tau_{b}<\tau_{a}\right)  $ and consider
column-vectors $\mathbf{h}_{k}\overset{\mathrm{def}}{=}(h_{k}(i))_{1\leq i\leq
m}$. In order to find $\mathbf{h}_{k}$ we introduce the matrices
$\varphi_{k+1}\overset{\mathrm{def}}{=}(\varphi_{k+1}(i,j))_{1\leq i,j\leq m}%
$, were
\begin{equation}
\varphi_{k+1}(i,j)\overset{\mathrm{def}}{=}Pr_{\omega,(k,i)}\left(  \tau
_{k+1}<\tau_{a},\ \xi(\tau_{k+1})=(k+1,j)\right)  . \label{DefPhi}%
\end{equation}
These matrices satisfy (\ref{phi1}) (with $a=0$) with the modified boundary
condition $\varphi_{a+1}=0$. The equation (\ref{EqD0.1}) with $\psi_{k}$'s
defined by (\ref{EqPsi}) now yields $\Delta_{k+1}=B_{k}...B_{a+1}\psi
_{a+1}\varphi_{a+2}...\varphi_{k+1}$ and hence
\begin{equation}
\left\Vert \Delta_{k+1}\right\Vert \leq\left\Vert B_{k}...B_{a}\right\Vert
\leq C\exp(\Phi_{k}-\Phi_{a}) \label{IneqDelta}%
\end{equation}
The Markov property also implies that $\mathbf{h}_{k}=\varphi_{k+1}%
\mathbf{h}_{k+1}$ and hence
\begin{equation}
\mathbf{h}_{k}=\varphi_{k+1}\varphi_{k+2}\ldots\varphi_{b}\mathbf{1}%
\ \hbox{ since }\ \mathbf{h}_{b}=\mathbf{1}. \label{hbyPhi}%
\end{equation}
We view the probabilities $Pr_{\left(  k,\cdot\right)  }\left(  \tau_{a}%
<\tau_{b}\right)  $ as the column vector $\mathbf{1}-\mathbf{h}_{k}$. Then,
presenting $\varphi_{b}=\psi_{b}-\Delta_{b}$, we can have%
\begin{align*}
Pr_{\left(  k,\cdot\right)  }\left(  \tau_{a}<\tau_{b}\right)   &
=\mathbf{1}-\varphi_{k}\ldots\varphi_{b-1}\mathbf{1}=\mathbf{1}-\varphi
_{k+1}\ldots\varphi_{b-1}(\psi_{b}-\Delta_{b})\mathbf{1}\\
&  =\mathbf{1}-\varphi_{k+1}\ldots\varphi_{b-1}\mathbf{1}+\varphi_{k+1}%
\ldots\varphi_{b-1}\Delta_{b}\mathbf{1}\\
&  \leq\mathbf{1}-\varphi_{k+1}\ldots\varphi_{b-1}\mathbf{1}%
+||\Delta_{b}||\mathbf{1.}%
\end{align*}
Iterating this inequality, we obtain that
\[
Pr_{\left(  k,\cdot\right)  }\left(  \tau_{a}<\tau_{b}\right)  \leq
\sum_{y=k+1}^{b}||\Delta_{y}||\mathbf{1}%
\]
and (\ref{hit1_1}) follows from (\ref{IneqDelta}). $\Box$

\begin{lemma}
\label{Le_Hittingtime_Exp}Let $a<b,$ and $\tau$ be the hitting time of
$L_{a}\cup L_{b}$ -- the union of two layers. Then if $a\leq k\leq b,$ we have%
\[
E_{\left(  k,i\right)  }\left(  \tau\right)  \leq C(b-a)^{2}\exp\left[
\min\left(  \sup_{a\leq s<t\leq b}\left(  \Phi\left(  s\right)  -\Phi\left(
t\right)  \right)  ,\sup_{a\leq s<t\leq b}\left(  \Phi\left(  t\right)
-\Phi\left(  s\right)  \right)  \right)  \right]
\]

\end{lemma}

\noindent\proof To prove that, consider column-vectors $\mathbf{e}_{k}=\left(
E_{(k,i)}\tau\right)  _{1\leq i\leq m}$. These vectors satisfy $\mathbf{e}%
_{a}=\mathbf{e}_{b}=\mathbf{0},$ and for $a<k<b:$%
\begin{equation}
\mathbf{e}_{k}=P_{k}\mathbf{e}_{k+1}+R_{k}\mathbf{e}_{k}+Q_{k}\mathbf{e}%
_{k-1}+\mathbf{1} \label{4.1}%
\end{equation}
To solve (\ref{4.1}), we use an induction procedure which allows us to find a
sequence of matrices $\varphi_{k}$ and vectors $\mathbf{d}_{k}$ such that
\begin{equation}
\mathbf{e}_{k}=\varphi_{k+1}\mathbf{e}_{k+1}+\mathbf{d}_{k}. \label{4.1.1}%
\end{equation}
Namely, we put $\varphi_{a+1}=0,\ \mathbf{d}_{a}=\mathbf{0}$ which according
to (\ref{4.1.1}) implies that $\mathbf{e}_{a}=\mathbf{0}$. Suppose next that
$\varphi_{k}$ and $\mathbf{d}_{k-1}$ are defined for some $k>a+1$. Then
substituting $\mathbf{e}_{k-1}=\varphi_{k}\mathbf{e}_{k}+\mathbf{d}_{k-1}$
into the main equation in (\ref{4.1}) we have
\[
\mathbf{e}_{k}=P_{k}\mathbf{e}_{k+1}+R_{k}\mathbf{e}_{k}+Q_{k}(\varphi
_{k}\mathbf{e}_{k}+\mathbf{d}_{k-1})+\mathbf{1}%
\]
and hence
\[
\mathbf{e}_{k}=(I-Q_{k}\varphi_{k}-R_{k})^{-1}\left(  P_{k}\mathbf{e}%
_{k+1}+Q_{k}\mathbf{d}_{k-1}+\mathbf{1}\right)
\]
which makes it natural to put
\begin{equation}
\varphi_{k+1}=(I-Q_{k}\varphi_{k}-R_{k})^{-1}P_{k} \label{4.19}%
\end{equation}
and
\begin{equation}
\mathbf{d}_{k}=B_{k}(\varphi_{k})\mathbf{d}_{k-1}+\mathbf{u}_{k},
\label{4.1.2}%
\end{equation}
where%
\[
\mathbf{u}_{k}=(I-Q_{k}\varphi_{k}-R_{k})^{-1}\mathbf{1},\ \ B_{k}(\varphi
_{k})=(I-Q_{k}\varphi_{k}-R_{k})^{-1}Q_{k}.
\]
The existence of matrices $\varphi_{k}$ follows from the fact that
$\varphi_{k}\geq0$ and $\varphi_{k}\mathbf{1}\leq\mathbf{1}$.

Iterating (\ref{4.1.1}) and (\ref{4.1.2}) we obtain
\[
\mathbf{e}_{k}=\mathbf{d}_{k}+\varphi_{k+1}\mathbf{d}_{k+1}+...+\varphi
_{k+1}...\varphi_{b-1}\mathbf{d}_{b-1}%
\]
and
\[
\mathbf{d}_{k}=\mathbf{u}_{k}+B_{k}(\varphi_{k})\mathbf{u}_{k-1}%
+...+B_{k}(\varphi_{k})...B_{a+1}(\varphi_{a+1})\mathbf{u}_{a}.
\]
Hence
\[
\left\Vert \mathbf{e}_{k}\right\Vert \leq\left\Vert \mathbf{d}_{k}\right\Vert
+\left\Vert \mathbf{d}_{k+1}\right\Vert +...+\left\Vert \mathbf{d}%
_{b-1}\right\Vert \leq C(b-k)\max_{k\leq j\leq b-1}||\mathbf{d}_{j}||.
\]
But $||B_{k}(\varphi_{k})...B_{l}(\varphi_{l})||\leq C\sup\nolimits_{a\leq
s<t\leq b}\exp\left(  \Phi\left(  s\right)  -\Phi\left(  t\right)  \right)  $
and therefore%
\[
E_{\left(  k,i\right)  }\left(  \tau\right)  \leq C(b-a)^{2}\exp\left[
\sup_{a\leq s<t\leq b}\left(  \Phi\left(  s\right)  -\Phi\left(  t\right)
\right)  \right]  .
\]

We obtain the same estimate with $\Phi$ replaced by $\Phi^{-},$ and using
(\ref{Compare_Phi_Phiprime}), we get the desired estimate. $\Box$

\begin{lemma}
\label{Le_Hittingtime2}Let $a\leq k<b$ and $\xi(t)$ be as in Lemma
\ref{Le_Hittingtime1}. Then for any $x>0$
\[
Pr_{\left(  k,i\right)  }\left(  \tau_{b}\geq x,\ \tau_{b}<\tau_{a}\right)
\leq\frac{C(b-a)^{2}}{x}\exp\left[  \sup\nolimits_{a\leq s<t\leq b}\left(
\Phi\left(  t\right)  -\Phi\left(  s\right)  \right)  \right]  .
\]

\end{lemma}

\noindent\proof Let again $\tau$ being the hitting time of $L_{a}\cup L_{b}$.
It is obvious that
\[
Pr_{\left(  k,i\right)  }\left(  \tau_{b}\geq x,\ \tau_{b}<\tau_{a}\right)
\leq Pr_{\left(  k,i\right)  }\left(  \tau\geq x\right)  .
\]
By the Markov inequality and Lemma \ref{Le_Hittingtime_Exp}, the result
follows. $\Box$

\begin{lemma}
\label{Le_Est_TransProb}Let $a<b,$ and consider the chain $\left\{  \xi
_{t}\right\}  $ on $\mathbb{S}_{a,b}$ with reflecting boundary conditions on
$a,b$, as above. Then for any $t\in\mathbb{N},$ $\left(  k,i\right)
,\ \left(  l,j\right)  \in\mathbb{S}_{a,b},$ we have%
\[
Pr_{\left(  k,i\right)  }\left(  \xi_{t}=\left(  l,j\right)  \right)  \leq
C\exp\left[  \Phi_{l}-\Phi_{k}\right]  .
\]

\end{lemma}

\noindent\proof%
\begin{align*}
\pi_{a,b}\left(  l,j\right)    & =\sum_{\left(  k^{\prime},i^{\prime}\right)
}\pi_{a,b}\left(  k^{\prime},i^{\prime}\right)  Pr_{\left(  k^{\prime
},i^{\prime}\right)  }\left(  \xi_{t}=\left(  l,j\right)  \right)  \\
& \geq\pi_{a,b}\left(  k,i\right)  Pr_{\left(  k,i\right)  }\left(  \xi
_{t}=\left(  l,j\right)  \right)
\end{align*}
for all $\left(  k,i\right)  ,$ $\left(  l,j\right)  \in\mathbb{S}_{a,b},$ and
all $t\in\mathbb{N}.$ The Lemma now follows with (\ref{Stat&Potential}).
$\Box$

We have now all the ingredients for the

\noindent\textbf{Proof of Proposition \ref{Prop_Main}}

We may assume that $0\in(a_{t},b_{t}].$ The case of $0\in\left(  b_{t}%
,c_{t}\right)  $ is handled similarly. We will write $a,b,c$ for $a_{t}%
,b_{t},c_{t},$ to simplify notations. We write $J_{t}$ for the interval
$\left[  b-\gamma\log^{2}t,b+\gamma\log^{2}t\right]  $

We have%
\begin{align}
Pr_{\left(  0,i\right)  }\left(  X\left(  t\right)  \notin J_{t}\right)   &
\leq Pr_{\left(  0,i\right)  }\left(  X\left(  t\right)  \notin J_{t},\tau
_{b}<\min\left(  \tau_{a},t\right)  \right)  +Pr_{\left(  0,i\right)  }\left(
\tau_{b}>\tau_{a}\right) \label{Main1}\\
&  +Pr_{\left(  0,i\right)  }\left(  \tau_{b}>t,\tau_{a}>\tau_{b}\right)
\nonumber
\end{align}

First we see that from Lemma \ref{Le_Hittingtime1}, and (\ref{Est1_Phi}),
(\ref{Est5_Phi}), (\ref{Est6_Phi})%
\begin{align}
Pr_{\left(  0,i\right)  }\left(  \tau_{b}>\tau_{a}\right)   &  \leq C\left(
b-a\right)  \exp\left[  \sup_{0\leq x\leq b}\Phi_{x}-\Phi_{a}\right]
\label{Main2}\\
&  \leq\frac{C\log^{2}t}{\sigma^{2}\delta}\exp\left[  -\delta\log t\right]
\leq t^{-\delta/2},\nonumber
\end{align}
if $t$ is large enough, and from Lemma \ref{Le_Hittingtime2} and
(\ref{Est3_Phi})%
\begin{align}
Pr_{\left(  0,i\right)  }\left(  \tau_{b}>t,\tau_{a}>\tau_{b}\right)   &
\leq\frac{C\log^{4}t}{t}\exp\left[  \sup\nolimits_{a\leq s<t\leq b}\left(
\Phi\left(  t\right)  -\Phi\left(  s\right)  \right)  \right] \label{Main3}\\
&  \leq\frac{C\log^{4}t}{t}\exp\left[  \left(  1-\delta\right)  \log t\right]
\leq t^{-\delta/2}.\nonumber
\end{align}
By the Markov property, we get%
\begin{equation}
Pr_{\left(  0,i\right)  }\left(  X\left(  t\right)  \notin J_{t},\tau_{b}%
<\min\left(  \tau_{a},t\right)  \right)  \leq\max_{s\leq t,1\leq j\leq
m}Pr_{\left(  b,j\right)  }\left(  X\left(  s\right)  \notin J_{t}\right)  .
\label{Main4}%
\end{equation}
Now%
\begin{equation}
Pr_{\left(  b,j\right)  }\left(  X\left(  s\right)  \notin J_{t}\right)  \leq
Pr_{\left(  b,j\right)  }\left(  \min\left(  \tau_{a},\tau_{c}\right)  \leq
t\right)  +Pr_{\left(  b,j\right)  }\left(  X^{\left(  a,c\right)  }\left(
s\right)  \notin J_{t}\right)  , \label{Main5}%
\end{equation}
where $X^{\left(  a,c\right)  }$ is the chain with reflecting boundary
conditions at $L_{a}$ and $L_{c}.$ The second summand is estimated by Lemma
\ref{Le_Est_TransProb} and (\ref{Est4_Phi}), which give%
\begin{equation}
Pr_{\left(  b,j\right)  }\left(  X^{\left(  a,c\right)  }\left(  s\right)
\notin J_{t}\right)  \leq C\exp\left[  \sup_{l\notin J_{t}}\Phi_{l}-\Phi
_{b}\right]  \leq Ct^{-\delta}\leq t^{-\delta/2}. \label{Main6}%
\end{equation}
To estimate the first summand in (\ref{Main5}) we observe that by
(\ref{Est5_Phi})%
\begin{align*}
Pr_{\left(  b-1,i\right)  }\left(  \tau_{a}<\tau_{b}\right)   &  \leq
C\exp\left[  -\Phi_{a}\right]  \left(  \exp\left[  \Phi_{b-1}\right]
+\exp\left[  \Phi_{b}\right]  \right)  \leq C\exp\left[  -\left(
1+\delta\right)  \log t\right] \\
&  \leq t^{-1-2\delta/3},
\end{align*}
and similarly%
\[
Pr_{\left(  b+1,i\right)  }\left(  \tau_{c}<\tau_{b}\right)  \leq
t^{-1-2\delta/3}.
\]
If, starting in $\left(  b,j\right)  ,$ the chain reaches $L_{a}$ or $L_{c}$
in time $t,$ there is at least one among the first $t/2$ of the excursions
from $L_{b}$ which reaches $L_{a}\cup L_{c}.$ By the above estimates, each
such excursion has at most probability $t^{-1-2\delta/3}$ to be
\textquotedblleft successful\textquotedblright, and therefore
\begin{equation}
Pr_{\left(  b,j\right)  }\left(  \min\left(
\tau_{a},\tau_{c}\right)  \leq t\right)  \leq1-\left(
1-t^{-1-2\delta/3}\right)  ^{t/2}\leq t^{-\delta/2}.
\label{Main7}%
\end{equation}
Combining (\ref{Main1})-(\ref{Main7}), we get%
\[
Pr_{\left(  0,i\right)  }\left(  X\left(  t\right)  \notin J_{t}\right)
\leq4t^{-\delta/2}.
\]
This proves the claim.

\section{Appendix}

Most (if not all) of the results in this Appendix are not new. The main reason
for including them is that we want to present them in the form which is needed
for our purpose; this is particularly relevant in the case of Markov chains
generated by contracting transformations. We also hope that a more
self-contained paper makes an easier reading.

\subsection{ The CLT and the invariance principle (IP) for stationary Markov
chains.}

We first recall, in subsection \ref{general-chains}, the classical results of
B. M. Brown \cite{Br} about the CLT and the IP for martingales. We then
explain in subsection \ref{general-chains1} that the reduction of the proof of
the CLT for Markov chains to the martingale case invented by Gordin and
Lifshits \cite{GL} can be easily extended to obtain the IP for Markov chains.
Finally, in subsection \ref{rand-tr}, we prove that the Gordin-Lifshits
conditions are satisfied for a class of Markov chains generated by contracting transformations.

\subsubsection{The CLT and the IP for martingales (by B. M. Brown \cite{Br}).
\label{general-chains}}

Let $\{\, S_{n},\ \mathcal{F}_{n}\,\},\ n=1,\ 2,...$ be a martingale on the
probability space $(\Omega,\mathcal{F},\mathbb{P})$. Put $U_{n}=S_{n}-S_{n-1}$
with $S_{0}=0$. The expectation with respect to $\mathbb{P}$ is denoted by
$\mathbb{E}$, and $\mathbb{E}_{j-1}$ stands for the conditional expectation
$\mathbb{E}(\cdot\,|\,\mathcal{F}_{j-1}) $. Let $\sigma_{n}^{2}=\mathbb{E}%
_{n-1}(U_{n}^{2})$, $V_{n}^{2}=\sum_{j=1}^{n}\sigma_{j}^{2}$, and $s_{n}%
^{2}=\mathbb{E}(V_{n}^{2})=\mathbb{E}(S_{n}^{2})$. The main assumption in
\cite{Br} concerned with martingales is:
\begin{equation}
\label{marting1}V_{n}^{2}s_{n}^{-2}\rightarrow1\ \ \hbox{ in probability as
$n\rightarrow\infty$}.
\end{equation}
We says that the Lindeberg condition holds for the class of martingales
satisfying (\ref{marting1}) if for any $\varepsilon>0$
\begin{equation}
\label{marting2}s_{n}^{-2}\sum_{j=1}^{n}\mathbb{E}U_{j}^{2}I(|U_{j}%
|\ge\varepsilon s_{n})\rightarrow0\ \ \hbox{ as \ $n\rightarrow\infty$},
\end{equation}
where $I(\cdot)$ is a characteristic function of a set.

For $t\in[\,0,1\,]$ define a sequence of piecewise linear random
functions
\begin{equation}%
\begin{array}
[c]{c}%
\label{defu} u_{n}(t)=s_{n}^{-1}\left(  S_{k}+U_{k+1}(ts_{n}^{2}-s_{k}%
^{2})(s_{k+1}^{2}-s_{k}^{2})^{-1}\right) \\
\hbox{ if $s_k^2\le ts_n^2\le s_{k+1}^2,\
k=0,1,...,n-1$.}
\end{array}
\end{equation}
The following two theorems from \cite{Br} describe the asymptotic behaviour of
the sequences $S_{n}$ and $u_{n}(\cdot)$.

\begin{theorem}
\label{CLT} If (\ref{marting1}) and (\ref{marting2}) hold, then $S_{n}$ is
asymptotically normal:
\begin{equation}
\label{marting3}\lim_{n\to\infty}\mathbb{P}\{s_{n}^{-1}S_{n}\le x\,\} =
(2\pi)^{-\frac{1}{2}}\int_{-\infty}^{x}e^{-\frac{1}{2}y^{2}}dy
\end{equation}
for all $x$. Furthermore, all finite dimensional distributions of $u_{n}(t)$
converge weakly, as $n\to\infty$, to those of a standard Wiener process $W(t)$
on $0\le t\le1$ (that is $W(0)=0$ and $\mathbb{E}W^{2}(1)=1$).
\end{theorem}

\begin{theorem}
\label{IP} Let $\{\, C[0,1],\mathcal{B},P_{W}\,\}$ be the probability space
where $C[0,1]$ is the space of continuous functions with the $\mathrm{sup}$
norm topology, $\mathcal{B}$ being the Borel $\sigma$-algebra generated by
open sets in $C[0,1]$, and $P_{W}$ the Wiener measure. Let $\{\mathbb{P}%
_{n}\}$ be the sequence of probability measures on $\{\, C[0,1],\mathcal{B}%
\,\}$ determined by the distribution of $\{\,u_{n}(t),\ 0\le t\le1\,\}$. Then
if (\ref{marting1}) and (\ref{marting2}) hold, $\mathbb{P}_{n}\rightarrow
P_{W}$ weakly as $n\to\infty$.
\end{theorem}

\subsubsection{The CLT and the IP for general Markov chains.
\label{general-chains1}}

In their famous work \cite{GL}, Gordin and Lifshits reduced the proof of the
CLT for Markov chains to that of martingales. They then applied the same
approach to the proof of the invariance principle for Markov chains in
\cite{GL1}. We shall explain their method here for the sake of completeness.

Let $z_{k}$, $k=1,2,...$, be a stationary ergodic Markov chain with a phase
space $(\mathfrak{X}, \mathcal{A})$, transition kernel $K(z,dy)$, and initial
distribution $\kappa$. Let $f:\mathfrak{X}\mapsto\mathbb{R}$ be a real valued
function on $\mathfrak{X}$ such that $\mathbb{E}f(z)=0$ and $\mathrm{Var}
f(z)<\infty$ (all expectations are taken with respect to the measure $\kappa
$). Let $L_{2}(\mathfrak{X}, \mathcal{A},\kappa)$ be the natural Hilbert space
associated with $\mathfrak{X}, \mathcal{A},\kappa$. By $\mathbf{I}$ we denote
the identity operator in this space, and by $\mathfrak{A}$ the transition
operator of the Markov chain: $\mathfrak{A}F(z)\overset{\mathrm{def}}{=}%
\int_{\mathfrak{X}}F(y)K(z,dy)$. Put
\begin{equation}
\label{defS}{S}_{n}=f(z_{1})+...+f(z_{n}) \hbox{ with the convention } {S}%
_{0}=0.
\end{equation}

\begin{theorem}
\label{clt-ip} Let $z_{k}$ be a Markov chain described above and suppose that
the function $f$ with $\mathbb{E}f=0$ can be presented as $f=(\mathbf{I}%
-\mathfrak{A})F$, where $F\in L_{2}(\mathfrak{X}, \mathcal{A},\kappa)$ and
$\mathbb{E}F=0$. Put $\sigma^{2}=||F||^{2}-||\mathfrak{A}F||^{2}%
\equiv\mathbb{E}F^{2} -\mathbb{E}(\mathfrak{A}F)^{2}$ and suppose that
$\sigma>0$. Then $\frac{{S}_{n}}{\sigma\sqrt n}$ converges in law towards the
standard Gaussian distribution $N(0,1)$ and the sequence ${S}_{n}$ satisfies
the invariance principle with parameter $\sigma$ in the sense of the
definition given in Section \ref{sec2.4}.
\end{theorem}

\proof Consider the identity which is due to Gordin (\cite{Gor}) and was used
by Gordin and Lifshits in \cite{GL}: $f(z_{k})=U(z_{k},z_{k+1})+F(z_{k}%
)-F(z_{k+1})$, where $U(z_{k},z_{k+1})=F(z_{k+1})-(\mathfrak{A}F)(z_{k})$.
This identity holds true because of the conditions imposed on $f$. Obviously,
$\mathbb{E}\{U(z_{k},z_{k+1})\,|\,z_{k}, ...,z_{1}\}=0$. Denote $U_{k+1}%
=U(z_{k},z_{k+1})$. In these notations we can write
\[
{S}_{n}=\hat{S}_{n}+F(z_{1})-F(z_{n+1}), \hbox{ where
$\hat{S}_n=\sum_{k=1}^{n}U_k$.}
\]
It is clear that if $\mathcal{F}_{n}$ is a $\sigma$-algebra generated by the
variables $z_{1},...,z_{n}$, then the sequence $\hat{S}_{n}$, $n=1,2,..$ is a
martingale with respect to the filtration $\mathcal{F}_{n}$, $n=1,2,...$. Let
us check that all conditions required by Theorems \ref{CLT} and \ref{IP} are
satisfied. Indeed, $\sigma_{j}^{2}=\mathbb{E}\{U_{j}^{2}\,|\,z_{j}\}
=(\mathfrak{A}F^{2})(z_{j})-[(\mathfrak{A}F)(z_{j})]^{2}$ is a stationary
sequence with $\mathbb{E}\sigma_{j}^{2}=||F||^{2}-||\mathfrak{A}F||^{2}%
=\sigma^{2}$. Relation (\ref{marting1}) takes the form
\[
(n\sigma^{2})^{-1}\sum_{j=1}^{n}\sigma_{j}^{2}\rightarrow1
\]
and is satisfied with probability 1 because of the Birkhoff Ergodic Theorem.
The Lindeberg condition (\ref{marting2}) takes the form
\[
\mathbb{E}U_{1}^{2}I(|U_{1}|\ge\varepsilon n\sigma^{2})\rightarrow
0\ \ \hbox{ as \ $n\rightarrow\infty$},
\]
and is obviously satisfied. Finally, functions (\ref{defu}) are now given by
\[
u_{n}(t)=n^{-\frac{1}{2}}\sigma^{-1}\left(  {S}_{k}+(tn-k)U_{k+1}\right)
\hbox{ if $k\le tn\le k+1,\
k=0,1,...,n-1$}
\]
and hence for $k\le tn\le k+1$
\[
v_{n}(t)=u_{n}(t)+n^{-\frac{1}{2}}\sigma^{-1}\left(  F(z_{1})-F(z_{k+1})
+(tn-k)(F(z_{k})-F(z_{k+1}))\right)  ,
\]
where $v_{n}(t)$ is as in (\ref{defv}). Since $F$ is square integrable and
$z_{n}$ is a stationary sequence, it follows that $n^{-\frac{1}{2}}\max_{1\le
k\le n}|F(z_{k})|\rightarrow0$ with probability 1 as $n\rightarrow\infty$.
Hence also the $\sup_{0\le t\le1}|v_{n}(t)-u_{n}(t)|\rightarrow0$ as
$n\rightarrow\infty$ with probability 1. All statements of our Theorem follow
now from Theorems \ref{CLT} and \ref{IP}. $\Box$

\subsubsection{The CLT and the IP for Markov chains
generated by contracting transformations. \label{rand-tr}}
Consider the following setup.

$(\Omega,\mathcal{F},\mathbb{P})$ is a probability space; the
related expectation is denoted $\mathbb{E}$.

$\mathrm{M}$ is a compact metric space equipped with a distance
$\rho(\cdot,\cdot)$.

$\mathfrak{B}$ is a semigroup of continuous Lipschitz
transformations of $\mathrm{M}$: for any $g\in \mathfrak{B}$ there
is a constant $l_g$ such that $\rho(g.y,g.y')\le l_g\rho(y,y')$
for any $y,\,y'\in\mathrm{M}$. Here and in the sequel $g.y$
denotes the result of the action of $g\in\mathfrak{B}$ on $y\in
\mathrm{M}$; this notation will be used most of the time but in
some cases we may write $g(y)$ rather than $g.y$.

For any $g_1,\,g_2\in\mathfrak{B}$ put
$\bar{\rho}(g_1,g_2)\de\sup_{y\in\mathrm{M}}\rho(g_1.y,g_2.y)$.
Obviously, $\bar{\rho}(\cdot,\cdot)$ defines a distance on
$\mathfrak{B}$. We can now consider a Borel sigma-algebra
generated by the corresponding open subsets of $\mathfrak{B}$;
this sigma-algebra will be denoted by $\mathfrak{S}$.

Consider a measurable mapping $g:\Omega\mapsto\mathfrak{B},\
\omega \mapsto g^\omega$ and for a $B\in \mathfrak{S}$ put
$\mu(B)\de \mathbb{P}\{\omega\,:\, g^\omega\in B\}$. We say that
$g$ is a random transformation of $M$. Let $g_k\in \mathfrak{B},\
k\ge1$ be a sequence if independent copies of $g$. Without loss of
generality we can assume that $g_k$ are defined on the same
probability space $(\Omega,\mathcal{F},\mathbb{P})$.

Denote by $\mathfrak{g}^{(j)}\de g_j\ldots g_1$ the product of
random transformations $g_1,...,g_j$ and let $\mu^{(j)}$ be the
probability distribution  of the product $\mathfrak{g}^{(j)}$.
This measure on $\mathfrak{B}$ is often called the $j^\mathrm{th}$
convolution power of the measure $\mu$ and is denoted by
$\mu^{(j)}=\mu^{\ast j}=\mu\ast\ldots\ast\mu$ ($j$ times).

A sequence of random transformations $g_k$ is said to be
\textit{contracting} if there are constants $C>0$ and $c,\ 0\le
c<1$ such that for any $y,\,y'\in\mathrm{M}$ and any $n\ge 1$
\begin{equation}\label{contr1}
\int_\mathfrak{B}\rho(g.y,g.y')\mu^{(n)}(dg)\equiv
\mathbb{E}\rho(g_{n}\ldots g_1.y ,g_{n}\ldots g_1.y')
 \le C c^n.
\end{equation}
\textit{Remark.} Perhaps it would be more natural to say that the
contraction property holds if
$\int_\mathfrak{B}\rho(g.y,g.y')\mu^{(n)}(dg) \le C c^n\rho(y,y')$.
However, (\ref{contr1}) is sufficient for our purposes and is what
we check in our applications.


\smallskip As usual, products of random transformations generate a
Markov chain with a state space $\mathrm{M}$. Namely, let
$\nu\equiv\nu(dy)$ be a probability measure on $\mathrm{M}$ and
let $y_1\in \mathrm{M}$ be chosen randomly according to the
distribution $\nu$ and independent of all $g_j$'s. For $k\ge1$
define $y_{k+1}\in\mathrm{M}$ by $y_{k+1}\de g_{k}.y_{k}\equiv
\mathfrak{g}^{(k)}.y_1$. The sequence of pairs $(g_k,y_k),\ k\ge1$
forms a Markov chain with a phase space
$\mathfrak{B}\times\mathrm{M}$; this chain will be denoted
$(\mathbf{g},\mathbf{y})$. Note that the $(\mathbf{y})$-component
of this chain, the sequences $y_k,\ k\ge1,$ is itself a Markov
chain with the phase space $\mathrm{M}$. Since $\mathrm{M}$ is a
compact space the chain $(\mathbf{y})$ has an invariant measure;
we shall suppose from now on that $\nu$ is such a measure which,
in turn, implies that $\mu(dg)\nu(dy)$ is an invariant measure of
the chain $(\mathbf{g},\mathbf{y})$. It is well known (and easy to
see) that if $g_k$ is a contracting sequence of random
transformations then the Markov chain $(\mathbf{y})$ has a unique
invariant measure.

Let $\mathcal{L}_2(\mathfrak{B}\times\mathrm{M})$ be the Hilbert
space of $\mu\times\nu$ square integrable real valued functions
and $\mathcal{C}(\mathfrak{B}\times\mathrm{M})$ be its subset of
continuous functions.

Given an $f\in \mathcal{C}(\mathfrak{B}\times\mathrm{M})$ let
${S}_n$ denote the related Birkhoff sums along a trajectory of the
Markov chain $(\mathbf{g},\mathbf{y})$:
\[{S}_n =
\sum_{k=1}^n f(g_k,y_k).
\]
By $\mathfrak{A}$ we denote the following Markov operator acting
in $\mathcal{L}_2(\mathfrak{B}\times\mathrm{M})$ and preserving
$\mathcal{C}(\mathfrak{B}\times\mathrm{M})$:
\begin{equation}\label{markov}
(\mathfrak{A}f)(g,y)\de
\int_{\mathfrak{B}}f(g^{\prime},g.y)\mu(dg^{\prime}).
\end{equation}
It follows from (\ref{markov}) that
\begin{equation}\label{mark}
(\mathfrak{A}^kf)(g,y)=
\int_{\mathfrak{B}\times\mathfrak{B}}f(g^{\prime},\tilde{g}g.y)
\mu(dg^{\prime})\mu^{(k-1)}(d\tilde{g}).
\end{equation}

\begin{theorem} \label{IPMC} Suppose that the sequence of random
transformations $g_k$ is contracting and $f$ is a continuous
bounded function on $\mathfrak{B}\times\mathrm{M}$ such that

(i) $\int_\mathfrak{B} f(g,y)\mu(dg)$ is Lipschitz on
$\mathrm{M}$, that is for some $C_f$
\[
|\int_\mathfrak{B} (f(g,y)-f(g,y'))\mu(dg)|\le C_f\rho(y,y'),
\]

(ii) $\int_\mathfrak{B} f(g,y)\mu(dg)\nu(dy) =0.$

\noindent Then the equation
\begin{equation}\label{markov1}
(I-\mathfrak{A})F= f,
\end{equation}
has a solution $F(g,y)$ which is continuous on $\mathfrak{B}\times
\mathrm{M}$ and
\[
\int_{\mathfrak{B}\times \mathrm{M}} F(g,y)\mu(dg)\nu(dy) =0.
\]
Besides, this solution is unique in
$\mathcal{L}_2(\mathfrak{B}\times\mathrm{M})$.

Denote by
\[
\sigma^2 =\int_{\mathfrak{B}\times \mathrm{M}}
(\mathfrak{A}F^2-(\mathfrak{A}F)^2)(g,y)\mu(dg)\nu(dy)
\]
If $\sigma >0$ then $\frac{{S}_n}{\sigma \sqrt n}$ converges in
law towards the standard Gaussian distribution $N(0,1)$ and the
sequence ${S}_n$ satisfies the invariance principle with parameter
$\sigma$.

\smallskip\noindent If $\sigma >0$ and, in addition to (i),
$|f(g,y)-f(g,y')|\le C_{f}(g)\rho(y,y')$ with $\int\log(1+
C_{f}(g))\mu(dg)<\infty$, then the invariance principle for the
sequence ${S}_n$ is satisfied uniformly in $y_1\in\mathrm{M}$.

\smallskip\noindent If $\sigma =0$, then the function $F(g,y)$
depends only on $y$ and for every $(g,y)$ in the support of
$\mu\times\nu$ one has
\begin{equation}\label{markov5} f(g,y) = F(y) -F(g.y).
\end{equation}
\end{theorem}
\proof \textit{The existence of $F$.} Equation (\ref{markov1}) can
be rewritten as $F=\mathfrak{A}F+f$ and, iterating this relation,
one obtains a formal series:
\begin{equation}\label{markov2}
F=\sum_{k=0}^{\infty} \mathfrak{A}^kf
\end{equation}
Condition (ii) of the Theorem and the invariance of the measure
$\mu(dg)\nu(dy)$ imply that
\[
\int_{\mathfrak{B}\times\mathrm{M}}(\mathfrak{A}^kf)(g,y)\mu(dg)\nu(dy)
=\int_{\mathfrak{B}\times\mathrm{M}}f(g,y)\mu(dg)\nu(dy)=0.
\]
Hence, the convergence in (\ref{markov2}) would follow if we prove
that
\begin{equation}\label{markov3}
|(\mathfrak{A}^kf)(g,y)- (\mathfrak{A}^kf)(\bar{g},\bar{y})| \le
\mathrm{const}\,c^{\frac{k}{n_0}}\ \hbox{for any $(g,y)$,
$(\bar{g},\bar{y})\in$  support of $\mu\times\nu$}.
\end{equation}
But it follows from (\ref{mark}) and condition (i) of the Theorem
that
\[
\begin{aligned}
&|(\mathfrak{A}^kf)(g,y)- (\mathfrak{A}^kf)(\bar{g},\bar{y})| \\
&  =\left|\int_{\mathfrak{B}}\left(\int_{\mathfrak{B}}
\left(f(g^{\prime},\tilde{g}g.y)-
f(g^{\prime},\tilde{g}\bar{g}.\bar{y})\right)
\mu(dg^{\prime})\right)\mu^{(k-1)}(d\tilde{g})\right|
\\ & \le C_f\int_{\mathfrak{B}} \rho(\tilde{g}g.y,
\tilde{g}\bar{g}.\bar{y}) \mu^{(k-1)}(d\tilde{g})\le C\,c^{n},
\end{aligned}
\]
where the last inequality is due to the contraction property
(\ref{contr1}). The existence and continuity of $F(g,y)$ is
proved.

\textit{Uniqueness.} As usual, to prove the uniqueness we have to
show that the homogeneous equation $F=\mathfrak{A}F$ has only a
trivial solution $F\equiv 0$ in the class of functions satisfying
the condition $\int_{\mathfrak{B}\times \mathrm{M}}
F(g,y)\mu(dg)\nu(dy) =0 $. To check that this is the case assume
that, to the contrary, there is an $F\in
\mathcal{L}_2(\mathfrak{B}\times\mathrm{M})$ such that
$F\not\equiv0$, satisfies the homogeneous equation, and has a zero
mean value. For a given $\epsilon>0$ find a  function $\tilde{F}$
which is Lipschitz on $\mathfrak{B}\times\mathrm{M}$ and
approximates $F$ in the sense that $||F-\tilde{F}||\le\epsilon$,
where $||\cdot||$ denotes the
$\mathcal{L}_2(\mathfrak{B}\times\mathrm{M})$ norm. The
$\tilde{F}$ can always be chosen so that $
\int_{\mathfrak{B}\times \mathrm{M}}
\tilde{F}(g,y)\mu(dg)\nu(dy)=0.$ Next, for any $n\ge 1$
\[
F=\mathfrak{A}^nF=\mathfrak{A}^n(F-\tilde{F})+\mathfrak{A}^n\tilde{F}.
\]
But then $\mathfrak{A}^n\tilde{F}\to 0$ uniformly in $(g,y)$ and
$||\mathfrak{A}^n(F-\tilde{F})||\le\epsilon$. Since $\epsilon$ can
be made arbitrarily small, we conclude that $F\equiv0$.

\textit{Proof of the CLT and the IP in the case $\sigma>0$.}
According to Theorem \ref{clt-ip} the existence of $F\in
\mathcal{L}_2(\mathfrak{B}\times\mathrm{M})$ satisfying equation
(\ref{markov1}) is the main condition under which both the Central
Limit Theorem and the Invariance Principle hold for Birkhoff sums
picked up along a realization of a trajectory of a Markov chain.
The ergodicity of the Markov chain is the other condition which is
needed and which in our case follows from the contraction
property. The CLT and the IP is thus proved.

\textit{Proof of the uniform IP in the case $\sigma>0$.} We write
${S}_n(y_1)$ for ${S}_n$ in order to emphasize the dependence of
this sequence on $y_1$. Clearly,
\begin{equation}\label{estS}
|{S}_n(y_1)-{S}_n(y_1')|\le \sum_{k=1}^n
|f(g_k,y_k)-f(g_k,y_k')|\le \sum_{k=1}^\infty
C_f(g_k)\rho(y_k,y_k').
\end{equation}
It follows from (\ref{contr1}) (due to the Chebyshev inequality)
that $\mathbb{P}$ almost surely $\rho(y_k,y_k')\le e^{-\varepsilon
k}$ for some $\varepsilon>0$ and $k\ge k(\varepsilon,\omega)$. It
is essential that $k(\varepsilon,\omega)$ does not depend on
$y_1,y_1'$. Next, due to the condition imposed on the function
$f$, the sequence $k^{-1}\log(1+C_f(g_k))\rightarrow0$ as
$k\rightarrow\infty$ $\mathbb{P}$ almost surely. Hence the right
hand side of (\ref{estS}) is $\mathbb{P}$ almost surely bounded
and the corresponding estimate does not depend on $y_1,y_1'$.

Let us now consider the dependence on $y_1$ of the relevant
$v_n(t)=v_n(t;y_1)$ (see (\ref{defv})). For $t\in[0,1]$, and $k\le
tn\le k+1,\ k=0,1,...,n-1$ we have:
\[
v_n(t;y_1)-v_n(t;y_1')=n^{-\frac{1}{2}}\left({S}_k(y_1)-S_k(y_1')+(f_{k+1}(y_1)-
f_{k+1}(y_1'))(tn-k)\right)
\]
with the obvious meaning of $f_{k+1}(y_1)$ and $f_{k+1}(y_1')$. It
is now clear that $\mathbb{P}$ almost surely
$v_n(t;y_1)-v_n(t;y_1')\rightarrow0$ as $n\rightarrow\infty$
uniformly in $y_1,y_1'$. This proves that the uniformity of the
invariance principle.

\textit{The case $\sigma=0$.} Note that
\[
(\mathfrak{A}F^2-\mathfrak{A}(F^2))(g,y)=\int_\mathfrak{B}\left(F(g',g.y)-
\int_\mathfrak{B}F(\tilde{g},g.y)\mu(d\tilde{g})\right)^2\mu(dg').
\]
Hence $\sigma=0$ implies that for $\mu\times\nu$-almost all
$(g,y)$ and $\mu$-almost all $g'$
\begin{equation}\label{markov4}
F(g',g.y)= \int_\mathfrak{B}F(\tilde{g},g.y)\mu(d\tilde{g}).
\end{equation}
But $F(\cdot,\cdot)$ is a continuous function of both variables
and hence (\ref{markov4}) holds for any $(g,y)$ from the support
of $\mu\times\nu$. This proves that $F$ depends only on the second
variable: $F(g',g.y)\equiv F(g.y)$ (we note that $g.y$ runs over
the whole of the support of $\nu$ when $(g,y)$ runs over the
support of $\mu\times\nu$). Finally, one obtains (\ref{markov5})
by substituting $F(y)$ (rather than $F(g,y)$) into
(\ref{markov1}). $\Box$

\subsubsection{Markov chains generated by contracting
transformations: characterization of the support of the invariant
measure.} The aim of this section is to give a characterization of
the support of an invariant measure of a Markov chain generated by
contracting transformations in terms of fixed points of these
transformations.

We work here within the same setup as in section \ref{rand-tr}.
This applies to the sequence $g_j$, $j\ge1$, the metric space
$(\mathrm{M},\rho)$, the semigroup $\mathfrak{B}$ of
transformations of $\mathrm{M}$, the Markov chain $y_j$ defined by
$y_{j+1}=g_j.y_j$, $j\ge 1$ (with $y_1$ being a random element
independent of all $g_j$'s). However, we shall suppose that
$\mathfrak{B}$ is generated by the transformations belonging to
the support $\mathcal{J}_0$ of the distribution $\mu$ of $g_j$'s.
This difference is important for Lemma \ref{support1}.

Let $\nu$ be the stationary measure of our chain and
$\mathrm{M}_0$ be the support of $\nu$.

As usual, we say that a transformation $g\in\mathfrak{B}$ is a
contraction on a subset $\mathrm{M}_0\subset\mathrm{M}$ if there
is an $n\ge1$ and a $c\in [0,1)$ (both $n$ and $c$ may depend on
$g$) such that $\rho(g^n.x',g^n.x'')\le c \rho(x',gx'')$ for any
$x'$, $x''\in \mathrm{M}_0$. If $g\in\mathfrak{B}$, then by $x_g$
we denote a fixed point of the transformation $g$: $g.x_g=x_g$.

\begin{lemma}\label{support0} If $g\in
\mathfrak{B}$ is a contraction on $\mathrm{M}$ then its fixed
point $x_g\in\mathrm{M}$, belongs to the support $\mathrm{M}_0$ of
the invariant measure $\nu$ of the Markov chain $y_j$.
\end{lemma}
\proof Consider a random infinite sequence $g_1,g_2,...$. Since
$g\in\mathcal{J}_0$, almost every such sequence has the property
that for any $k\ge1$ and any $\delta>0$ there are infinitely many
$i$'s such that each element of the part $g_{i},...,g_{i+nk-1}$ of
the sequence approximates $g$ so closely that
\[
\bar{\rho}(g^{nk},\mathfrak{g}_{i}^{(nk)})\le \delta\ \hbox{ where
}\ \mathfrak{g}_{i}^{(nk)} \de g_{i+nk-1}...g_{i}.
\]
Moreover, by the law of large numbers these $i$'s have a positive
frequency. Since
\[
\rho(x_g,g^{nk}.x')= \rho(g^{nk}x_g,g^{nk}.x') \le c^k\rho(x_g,x')
\]
for any $x'\in \mathrm{M}$, we have that
\[
\rho(x_g,\mathfrak{g}_{i}^{(nk)}.x')\le c^k\rho(x_g,x')+
\rho(g^{nk}.x',\mathfrak{g}_{i}^{(nk)}.x')\le
c^k\rho(x_g,x')+\delta.
\]
Hence any (small) neighbourhood of $x_g$ is visited by the
sequence $\mathfrak{g}_{1}^{(j)}.x'$, $j\ge1$, infinitely many
times and, moreover, this happens with a positive frequency for
almost every sequence $g_j$, $j\ge1$. This implies that $x_g\in
\mathrm{M}_0$ and $(g,x_g)\in\mathcal{J}_0\times \mathrm{M}_0$.
$\Box$

Note that if the invariant measure $\nu$ of our Markov chain is
ergodic, then the support $\mathrm{M}_0$ of this measure is a
minimal set of $\mathfrak{B}$. The latter by definition means that
the orbit $\{g.x:\, g\in\mathfrak{B}\}$ of any $x\in\mathrm{M}_0$
is everywhere dense in $\mathrm{M}_0$.
\begin{lemma}\label{support1} Let $\mathrm{M}_0\subset\mathrm{M}$
be a minimal set of $\mathfrak{B}$. Suppose that there exist a
$\hat{g}\in\mathfrak{B}$ which is a contraction on $\mathrm{M}_0$.
Consider the set of all fixed points of $\mathfrak{B}$ belonging
to $\mathrm{M}_0$:
\[
\mathrm{Fix}_{\mathrm{M}_0}(\mathfrak{B})\de \{x:\,
x\in\mathrm{M}_0 \hbox{ and there is a $g \in \mathfrak{B}$ such
that g.x=x } \}.
\]
Then $\mathrm{Fix}_{\mathrm{M}_0}(\mathfrak{B})$ is everywhere
dense in $\mathrm{M}_0$.
\end{lemma}
\proof The contraction $\hat{g}$ given to us by the condition of
the Lemma has a fixed point $\hat{x}\in \mathrm{M}_0$ (it may have
other fixed points too, but we are interested only in this one).
Since $\mathrm{M}_0$ is minimal it coincides with the closure of
the orbit $\{g.\hat{x}:\, g\in\mathfrak{B}\}$. For a given
$g\in\mathfrak{B}$ let us consider the point $g.\hat{x}$. We shall
now show that for a sufficiently large $n$ the transformation
$g\hat{g}^{n}$ has a fixed point which we shall denote
$x_{g\hat{g}^{n}} $. Indeed, for any $x',x''\in\mathrm{M}_0$
\[ \rho(g\hat{g}^{n}.x',g\hat{g}^{n}.x'')\le
l_g\rho(\hat{g}^{n}.x',\hat{g}^{n}.x'')\le l_gc^n\rho(x',x'').
\]
If $n$ is such that $l_gc^n<1$, then there is a fixed point
$x_{g\hat{g}^{n}}$ of $g\hat{g}^{n}$. On the other hand, it is
obvious that $g\hat{g}^{n}.x'\rightarrow g.\hat{x}$ as
$n\to\infty$ uniformly in $x'\in\mathrm{M}_0$ because
$\hat{g}^{n}.x'\rightarrow \hat{x}$ uniformly in
$x'\in\mathrm{M}_0$. It follows that in particular
$x_{g\hat{g}^{n}}\rightarrow g.\hat{x}$ and this proves the Lemma.
$\Box$

\subsection{Products of positive matrices.} Lemma \ref{pos} below
explains two versions of a well known contraction property of
products of positive matrices (see, e.g. \cite{FK}). The first
version of this property has already been explained and proved in
the Appendix to \cite{BG} and we therefore prove here only the
second version. There is a slight difference in the notations used
in this paper and those we have introduced in \cite{BG} and no
difference in the proof; we emphasize once again that this is done
for the purposes of completeness and convenience of references in
the proofs of other theorems.
\begin{lemma}\label{pos}
Let $a_{n}=(a_{n}(i,j)),\ n=1,2,\ldots$ be a sequence of positive
$m\times m$ matrices, $a_{n}>0$. Put $ \tilde{H}_{n}\de
a_{n}a_{n-1}\ldots a_{1}$, $ H_{n}\de a_{1}a_{2}\ldots a_{n}$ and
denote
\[
\tilde{\delta}_{r}=\min_{i,j,k}a_{r}(i,j)a_{r-1}(j,k)(\sum_{j}a_{r}(i,j)a_{r-1}
(j,k))^{-1},\ \ 2\le r\le n\]
\begin{equation}\label{delta}
{\delta}_{r}=\min_{i,j,k}a_{r}(i,j)a_{r+1}(j,k)(\sum_{j}a_{r}(i,j)
a_{r+1}(j,k))^{-1},\ \ 1\le r\le n-1.
\end{equation}
Suppose that
\[
\sum_{r=2}^{\infty}\tilde{\delta}_{r}=\infty
\]
Then the products $H_{n}$ and $\tilde{H}_{n}$ can be presented as
follows:
\begin{equation}\label{A1}
H_{n}=D_{n}[\left( c_n(1)\mathbf{1},\ldots,c_n(m)\mathbf{1}\right)
+\phi_{n}], \ \ \ \tilde{H}_{n}=\tilde{D}_{n}[\left(
\tilde{c}(1)\mathbf{1},\ldots,\tilde{c}(m)\mathbf{1}\right)
+\tilde{\phi}_{n}],
\end{equation}
where:

$D_{n}$ and $\tilde{D}_{n}$ are diagonal matrices with positive
diagonal elements;

$\left\| \phi_{n}\right\| \leq\prod_{r=1}^{n-1}(1-m\delta_{r})$,
$\left\| \tilde{\phi}_{n}\right\|
\leq\prod_{r=2}^{n}(1-m\tilde{\delta}_{r})$;

$\tilde{c}(j)$ are strictly positive numbers which are uniquely
defined by the sequence $\{a_{k}\}_{k\ge1}$, do not depend on $n$,
and such that $\sum_{j}\tilde{c}(j)=1$;

$c_n(j)$ are strictly positive numbers with $\sum_{j}c_n(j)=1$
(note that $c_n(j)$, unlike the $\tilde{c}(j)$,  do depend on $n$
and, generally, do not have a limit).
\end{lemma}

\proof Present $H_{n}$ as follows:
\[
H_{n}=D_{n}D_{n}^{-1}a_{1} D_{n-1}D_{n-1}^{-1}a_{2}\ldots
D_{1}^{-1}
a_{n}=D_{n}\tilde{a}_{1}\tilde{a}_{2}\ldots\tilde{a}_{n}%
\]
where  $\tilde{a}_{r}\equiv D_{n-r+1}^{-1}a_{r}D_{n-r}$,
$\tilde{D}_{0}\de I$, and
$D_{n-r}=\mathrm{diag}\left(D_{n-r}(1),...,D_{n-r}(m)\right)$ are
diagonal matrices, with $D_{n-r}(i)$ chosen so that to make
matrices $\tilde{a}_{r}$ stochastic. It is very easy to see that
the only such choice is given by
\[
D_{n-r}(i)=\sum_{i_{r+1},\ldots,i_{n}}a_{r+1}\left(i,i_{r+1}\right)a_{r+2}\left(
i_{r+1},i_{r+2}\right)  \ldots a_{n}\left(  i_{n-1},i_{n}\right)
\]
and
\begin{align}
\tilde{a}_{r}(i,j)
=\frac{a_{r}(i,j)\sum_{i_{r+1},\ldots,i_{n}}a_{r+1}\left(j,i_{r+1}\right)
\ldots a_{n}\left(
i_{n-1},i_{n}\right)}{\sum_{i_r,i_{r+1},\ldots,i_{n}}a_{r}\left(i,i_{r}\right)a_{r+1}\left(
i_{r},i_{r+1}\right)  \ldots a_{n}\left(
i_{n-1},i_{n}\right)}\geq\delta_{r}. \label{delta1}
\end{align}
It is well known that the last estimate implies the following
presentation of the product of stochastic matrices
$\tilde{a}_{n}$:
\[
\tilde{a}_{1}\tilde{a}_{2}\ldots\tilde{a}_{n}=(c_n(1)\mathbf{1},\ldots
,c_n(m)\mathbf{1})+\phi_{n},%
\]
where
\begin{equation}\label{delta2}
\min_i\tilde{a}_{n}(i,j)\le c_n(j)\le\max_i\tilde{a}_{n}(i,j)
\end{equation}
and the matrices $\phi_{n}$ are such that
\[
\left\| \phi_{n}\right\| \leq\prod_{r=1}^{n-1}(1-m\delta_{r}).
\]
$\Box$

\subsection{A stability estimate.} The stability property
which we explain below is definitely well known to specialists in
the relevant field. Given that the proof is very short, it seems
that it is easier for us to prove it than to find a relevant
reference.

Let $b_n$ and $b_n'$ be two sequences of transformations of a
metric space $(\mathbb{X},\mathfrak{r})$ and $x_{n+1}\de
b_n(x_n)$, $x_{n+1}'\de b_n'(x_n')$, $n\ge1$,  with given initial
values $x_1,\ x_1'\in \mathbb{X}$. For any two transformations $b$
and $b'$ put
$\bar{\rho}(b,b')\de\sup_{x\in\mathbb{X}}\mathfrak{r}(b(x),b'(x))$

\begin{lemma}\label{stability}
Suppose that

(a) $b_n$ are uniformly contracting, that is there is a $c,\ 0\le
c<1,$ such that for any $x,y\in\mathbb{X}$ we have
$\mathfrak{r}(b_n(x),b_n(y))\le c \mathfrak{r}(x,y)$;

(b) $\bar{\rho}(b_n,b_n')\rightarrow0$ as $n\rightarrow\infty$.

Then $\mathfrak{r}(x_n,x_n')\rightarrow0$ as $n\rightarrow\infty$.

If, instead of (b), a stronger property holds, namely
$\bar{\rho}(b_n,b_n')\le C_2 c_0^n\bar{\rho}(b_1,b_1')$ for some
$C_2$ and $c_0<1$, then for  $\epsilon>0$ there is a constant
$C_3$ such that
\begin{equation}\label{estr}
\mathfrak{r}(x_n,x_n')\le C_3
\tilde{c}^n(\bar{\rho}(b_1,b_1')+\mathfrak{r}(x_1,x_1')), \hbox{
where } \tilde{c}=\max(c,c_0)+\epsilon.
\end{equation}
\end{lemma}
\proof Put $d_n\de \bar{\rho}(b_n,b_n')$ and
$r_n\de\mathfrak{r}(x_n,x_n')$. Since
\[
\begin{aligned}
\mathfrak{r}(x_{n+1},x_{n+1}')=&\mathfrak{r}(b_n(x_n),b_n'(x_n'))\le
\mathfrak{r}(b_n(x_n),b_n(x_n'))+\mathfrak{r}(b_n(x_n'),b_n'(x_n'))\\
&\le c\mathfrak{r}(x_n,x_n')+\bar{\rho}(b_n,b_n'),
\end{aligned}
\]
we have that
\begin{equation}\label{Eqr}
\begin{aligned}
r_{n+1}\le cr_n+d_n\le d_n+cd_{n-1}+...+c^kd_{n-k}
+c^{k+1}r_{n-k}.
\end{aligned}
\end{equation}
For a given $\epsilon>0$ choose $k$ so that $c^kr_{n-k}\le
\epsilon$ (which is possible because $\mathbb{X}$ is a compact
space and thus $r_{n-k}$ is a uniformly bounded sequence). Next
choose $N(\epsilon, k)$ so that $d_{n-j}\le \epsilon$ when $n-j\ge
N(\epsilon, k)-k$. It follows now from (\ref{Eqr}) that
$r_{n}\le(2-c) (1-c)^{-1}\epsilon$ when $n>N(\epsilon, k)$. This
proves the first statement of the Lemma.

To prove the second statement substitute $k=n$ into (\ref{Eqr})
and take into account the stronger estimates for $d_n$. Estimate
(\ref{estr}) follows with an evident choice of $C_3$. $\Box$

\textit{Remark.} The second statement of this Lemma does not use
the fact that $\mathbb{X}$ is a compact space.

\sni\textbf{Acknowledgements.} This work was supported by the
following grants of the Swiss National Foundation: 200020-107739/1
and 200020-116348. We are grateful to the Isaac Newton Institute for
its hospitality during the program {\it Interaction and Growth in
Complex Stochastic Systems} held in Cambridge, UK in 2003. We also
thank the European Science Foundation Research Networking Programme
on {\it Phase-Transitions and Fluctuation Phenomena for Random
Dynamics in Spatially Extended Systems (RDSES)} for its financial
support.

\bigskip

Erwin Bolthausen, Universit\"{a}t Z\"{u}rich, Institut f\"{u}r
Mathematik, Winterthurerstrasse 190, CH-8057 Z\"{u}rich

email: eb@math.unizh.ch

\bigskip

Ilya Ya. Goldsheid, School of Mathematical Sciences, Queen Mary
and Westfield College, University of London, London E1 4NS, UK

email: I.Goldsheid@qmul.ac.uk
\end{document}